\documentclass{article}
\pdfoutput=1

\usepackage[utf8]{inputenc}
\usepackage[T1]{fontenc}
\usepackage[english]{babel}
\usepackage{amsthm}
\usepackage{amsmath}
\usepackage{amsfonts}
\usepackage{amssymb}
\usepackage{graphicx}
\usepackage{mathtools}
\usepackage{hyperref}
\usepackage[usenames,dvipsnames]{xcolor}

\usepackage[shortlabels]{enumitem}

\title{Minimal degree rational open up mappings 
and related questions}

\author{Sergei Kalmykov \and Béla Nagy \and Olivier S\`ete}

\date{}

\usepackage{amsopn}

\newcommand*{\C}{\mathbf{C}}

\newcommand*{\bC}{\mathbf{C}}

\DeclareMathOperator{\interior}{int}
\DeclareMathOperator{\exterior}{ext}

\DeclareMathOperator{\coeff}{coeff}
\DeclareMathOperator{\id}{id}

\DeclarePairedDelimiter{\abs}{\lvert}{\rvert}
\DeclarePairedDelimiter{\cc}{[}{]}
\newcommand{\coloneq}{\mathrel{\mathop:}=}
\newcommand{\cp}{\zeta}
\newcommand{\cv}{\eta}

\theoremstyle{plain}
\newtheorem{theorem}{Theorem}[section]
\newtheorem{lemma}[theorem]{Lemma}
\newtheorem{proposition}[theorem]{Proposition}
\newtheorem{corollary}[theorem]{Corollary}

\theoremstyle{definition}
\newtheorem*{conjecture}{Conjecture}

\newtheorem{remark}[theorem]{Remark}

\begin{document}
\maketitle

\begin{center}
Dedicated to Vilmos Totik on his 70th birthday
\end{center}

\begin{abstract}
We establish the existence and uniqueness of rational conformal maps of minimal 
degree $n+1$ for opening up $n$ arcs.
In earlier results, the degree was exponential in $n$.
We also discuss two related problems.
(a) We establish existence of rational functions of minimal degree with 
prescribed critical values, and show that the number of (suitably normalized) 
rational functions is given in terms of the Hurwitz numbers.
(b) We consider the problem of finding rational functions of minimal degree 
with prescribed critical points, where we establish existence of solutions by 
considering certain polynomial equations, and where the number of normalized 
solutions is bounded from above by a Catalan number.
We illustrate our results with two examples.


MSC 2020 Classification codes:
30C20, 30C10, 13P15

\end{abstract}

\section{Introduction}
\label{sect:introduction}

We call a rational function $F$ of \emph{type} $(n,m)$ if it can be written as
$F = P/Q$, where $P$ and $Q$ are coprime polynomials with $\deg(P) = n$
and $\deg(Q) = m$.  The \emph{degree} of $F$ is $\deg(F) = \max \{ n, m \}$.
We denote the extended complex plane by $\bC_\infty = \bC \cup \{ \infty \}$.

The original purpose of this research was to prove the following theorem.

\begin{theorem} \label{thm:openup}
Let $\gamma_1, \ldots, \gamma_n$ be disjoint 
Jordan arcs in the complex plane.
Then there exists a rational function $F$ of type $(n+1,n)$ and a compact set 
$K \subset \bC$ bounded by $n$ disjoint Jordan curves 
such that $F$ is a conformal map from $\bC_\infty \setminus K$
onto $\bC_\infty \setminus \cup_{j=1}^n \gamma_j$ and $F(\infty)=\infty$.
Moreover, $F$ and $K$ are unique up to pre-composition of $F$ with a linear 
transformation.
In particular, the normalization $F(z) = z + O(1/z)$
at infinity determines $F$ and $K$ uniquely.
\end{theorem}

Recall that Jordan arcs and Jordan curves are homeomorphic images of the
interval $[0,1]$ and the unit circle, respectively.
Throughout this article, we always mean that a conformal map is also injective.
Note that $F(\infty) = \infty$ always holds if $F$ is of type $(n+1, n)$.
The inverse function $F^{-1} : \bC_\infty \setminus \cup_{j=1}^n \gamma_j \to
\bC_\infty \setminus K$ ``opens up'' the arcs $\gamma_j$.  
We call $F^{-1}$, and for simpler language also  $F$, an \emph{open up mapping}
for the arcs $\gamma_1, \ldots, \gamma_n$.

Such results can be applied, among others, to
prove asymptotically sharp Bernstein- and Markov-type
inequalities on several Jordan arcs.
In~\cite{KalmykovNagy2015} and~\cite{KalmykovNagyTotik2017} it was shown how it works in the case of one arc.

For a rational function $F$ in Theorem~\ref{thm:openup},
the endpoints
of the arcs are critical values (see Proposition~\ref{prop:endpts_are_critvals}
below).
This observation leads to the problem of finding rational functions 
of minimal degree with prescribed critical values.
As for terminology, the \emph{critical points}
are the points in the set $\{ z \in \bC : F'(z) = 0 \}$, and
the \emph{critical values} are the elements of
$\{ F(z) : z \in \bC \text{ with } F'(z) = 0 \}$.

\begin{theorem} \label{thm:critvals}
Let $\cv_1,\cv_2,\ldots,\cv_{2n} \in \bC$ be distinct.
Then there exists a rational function $F$ of type $(n+1,n)$ such that the
set of critical values of $F$ is $\{ \cv_1,\cv_2,\ldots,\cv_{2n} \}$.
Moreover, each function can be normalized by $F(z) = z + O(1/z)$
at infinity.
If $n = 1$, there is exactly one normalized function.  If $n \geq 2$,
the number of normalized functions is $(n+1) H_n$ with the Hurwitz numbers
\begin{equation} \label{eqn:hurwitz_numbers}
H_n = 
\begin{cases}
\frac{(2n)!}{n!} (n+1)^{n-3}, & \text{if } n \geq 2, \\
1, &\text{if } n=1.
\end{cases}
\end{equation}
\end{theorem}

In contrast to Theorem~\ref{thm:openup},
the normalization $F(z) = z + O(1/z)$ at infinity does not uniquely determine
a rational function with prescribed critical values when $n \geq 2$.
We give an example for this in Section~\ref{sect:example}.

For each critical value $\cv_j$ of $F$, there exists a critical point $\cp_j$
with $F(\cp_j) = \cv_j$ and $F'(\cp_j) = 0$.  This remark 
leads to the related problem of finding a rational function 
of minimal degree with prescribed critical points. 
It turns out that this is a simpler problem 
(with half as many equations and unknowns) 
than the previous one and it is answered by the following theorem.

\begin{theorem} \label{thm:critpts}
Let $\cp_1, \ldots, \cp_{2n} \in \bC$ be distinct.
Then there exists a rational function 
$F$ of type $(n+1, n)$ such that
the set of critical points 
of $F$ is $\{ \cp_1, \ldots, \cp_{2n} \}$.
Moreover, each function $F$ can be normalized by
$F(z) = z + O(1/z)$
at infinity, and the number of normalized functions is bounded 
from above
by
the Catalan number $C_n = \frac{1}{n+1} \binom{2n}{n}$.
\end{theorem}

That the degree is, indeed,
minimal in Theorems~\ref{thm:openup}--\ref{thm:critpts} is shown next.

\begin{proposition} \label{prop:minimal}
The rational functions in Theorems~\ref{thm:openup},
\ref{thm:critvals} and~\ref{thm:critpts} are minimal
in the sense that if $F = P/Q$
with coprime polynomials $P, Q$,
then neither $\deg(P) < n+1$ nor $\deg(Q) < n$ can occur.
Moreover, they have only simple poles.
\end{proposition}

\begin{proof}
Since a rational function in Theorem~\ref{thm:openup} or~\ref{thm:critvals}
has $2n$ distinct critical points, it is enough to show the proposition
for rational functions in Theorem~\ref{thm:critpts} (i.e., with $2n$ distinct critical points).
Let $F$ be a rational function as in Theorem~\ref{thm:critpts}.
Assume that $F = P/Q$ with coprime polynomials $P$, $Q$, and that 
\begin{enumerate}
\item $\deg(F) \leq n$, or
\item $\deg(P) \leq n+1$ and $\deg(Q) \leq n$, and we have strict inequality for $P$ or $Q$.
\end{enumerate}
Then $F' = (P'Q-PQ')/Q^2$ and $\deg(P'Q-PQ') < 2n$.
In Theorem~\ref{thm:critpts} if $F$ has $2n$ distinct critical points then
we obtain $F' \equiv 0$, i.e., $F$ is constant, which is impossible.

It remains to show that $F$ has simple poles.
If $F = P/Q$ with $\deg(P) = n+1$ and $\deg(Q) = n$ has a pole
of order $m \geq 2$,
then, after cancelling common factors, the numerator of $F' = (P'Q - PQ')/Q^2$,
has degree at most $2n+1-m < 2n$ but $2n$ zeros, which is impossible.
\end{proof}

\begin{remark} \label{rem:lineartransformation}
Each function in Theorems~\ref{thm:openup}, \ref{thm:critvals} or~\ref{thm:critpts} (without normalization)
yields infinitely many solutions of the respective problem by suitable affine
transformations.
Indeed, let $a, b \in \bC$ with $a \neq 0$.
\begin{enumerate}
\item In Theorems~\ref{thm:openup} and~\ref{thm:critvals},
if $F$ is a solution, then so is $F(az+b)$.
\item In Theorem~\ref{thm:critpts}, if $F$ is a solution, then
$a F(z)+b$ is also a solution.
\end{enumerate}
In particular, if a solution exists, it is not unique.
We call two solutions $F, G$ equivalent, if there exist $a, b \in \bC$ with
$a \neq 0$ such that $G(z) = F(az+b)$ (in Theorems~\ref{thm:openup}
and~\ref{thm:critvals}) or $G(z) = a F(z) + b$ (in Theorem~\ref{thm:critpts}).
We can specify any normalization at infinity to obtain a unique
representative of each equivalence class:
Given $\alpha, \beta \in \bC$ with $\alpha \neq 0$ and a solution $F$ of
Theorem~\ref{thm:openup}, there exist unique $a, b$ as above, such that
$F(az+b) = \alpha z + \beta + O(1/z)$ for $z \to \infty$.
A similar statement holds for the other two theorems.
\end{remark}

The paper is organized as follows.
We give an overview of known related results in Section~\ref{sect:overview}.
In Section~\ref{sect:openup}, we show the existence and
uniqueness in Theorem~\ref{thm:openup} using the theory of Riemann surfaces.
This approach is purely geometric: it is short and intuitive,
but not constructive.
The results in Theorems~\ref{thm:critvals} and~\ref{thm:critpts} are
of algebraic nature, and we reformulate both in terms of polynomial systems
of equations.
In Section~\ref{sect:critvals}, we prove Theorem~\ref{thm:critvals}.
We give a geometric existence proof using the open up mapping, while the proof 
for the number of solutions builds on earlier results by Hurwitz and 
Mednykh on the number of Riemann surfaces with simple branch points.
In Section~\ref{sect:critpts}, we prove Theorem~\ref{thm:critpts} and use
algebraic tools to show the existence of solutions.
We conclude with a discussion in Section~\ref{sect:discussion}
and two examples in Section~\ref{sect:example}.

\section{Overview of some known, earlier results}
\label{sect:overview}

There are several results related to the three theorems above
that are scattered through the literature.
They appeared in various fields of mathematics and occurred 
in almost every decade in the last century.
Let us recall some of them,
not necessarily in chronological order. 

Starting with Theorem~\ref{thm:openup}, Widom's seminal paper~\cite{Widom}
must be mentioned where he 
iterated Joukowskii mappings to
construct a rational open up mapping,
see~\cite[pp.~206--207]{Widom}.
In Widom's construction, the rational function has degree $2^n$, i.e.,
it grows exponentially with the number of arcs.
Later, Widom's iterated construction appeared in connection with Riemann
surfaces, see the papers by Seppälä~\cite{Seppala}
and Seppälä and Hidalgo~\cite{HidalgoSeppala}.
Seppälä attributes this approach to Myrberg~\cite{Myrberg}, 
who, in turn, credits this idea to Poincaré, see~\cite[p.~4]{Myrberg}.

The question about the existence of rational functions 
with prescribed critical values
can be considered in general: is it possible
to cover Riemann surfaces with prescribed ramification sets?
For results in this direction and going back to a problem of Hurwitz,
we refer to 
Mednykh's paper~\cite{Mednykh}
and the references therein.

These ramification sets or branching points naturally lead to
Theorem~\ref{thm:critvals}.
Instead of rational functions, polynomials with prescribed critical values
were also investigated in Thom's paper~\cite{Thom}, in which
the existence of such polynomials was established.  See also
the papers by Mycielski and Paszkowski~\cite{MycielskiPaszkowski},
Kammerer~\cite{Kammerer1961}, Kuhn~\cite{Kuhn1969}
and Kristiansen~\cite{Kristiansen}
for the real case and further references.
Let us remark that Beardon, Carne, and Ng in~\cite{BCN}
investigated the properties of the mapping from critical points 
to critical values, realizing and describing 
a natural connection between the two problems in the class of polynomials.
This leads to Theorem~\ref{thm:critpts}.

Theorem~\ref{thm:critpts} was proved by Goldberg 
in~\cite{Goldberg} using projective spaces,
Grassmann manifolds, and homology classes.
She also counted the number of solutions.
Since we are also interested in
obtaining the solutions, we show the existence
by more constructive means.

For real rational functions with prescribed real critical points 
we refer to the paper of Eremenko and Gabrielov \cite{EremenkoGabrielov}
and the recent article~\cite{PeltolaWang}
where such rational functions
are used in 
Schramm--Loewner 
evolution.

An interesting application of the open up mapping is the computation of the logarithmic capacity of a compact set $E$ consisting of $n$ disjoint arcs.  
The rational function maps the exterior of $E$ to a domain with smooth
boundary, from which the logarithmic capacity of $E$ can be computed numerically with
a conformal map of Walsh, as described by Nasser, Liesen and S\`ete 
in~\cite{NasserLiesenSete} and~\cite{LiesenSeteNasser2017}.

The necessity of obtaining conformal representations 
by rational functions also appeared in the study of
multiple orthogonal polynomials, see~\cite{AKLT,ATvA,ATY,Lagomasino}.
In general, such representations are 
different from open up mappings,
but they are the same in the case of two arcs.
The case of two real intervals was considered 
by L\'opez Lagomasino, Pestana, Rodr\'{\i}guez, and Yakubovich 
in~\cite{Lagomasino}.

Instead of (general) rational functions,
similar questions can be considered
among (finite) Blaschke products. 
See, e.g.,
\cite{KrausRoth2008}, \cite{SemmlerWegert} 
for further references.
A similar open problem (determining a Blaschke
product or its zeros from critical values)
is also of interest and is raised in~\cite{SemmlerWegert}.

We mention that Theorem~\ref{thm:openup} is more precise
than~\cite[Prop.~5]{KalmykovNagy2015}.
Note that there is a minor flaw in the proof of existence
in~\cite[Prop.~5]{KalmykovNagy2015}, while the remaining proof is correct.
In our paper, we prove the existence of $F$ and $K$ with a different method.

\section{Existence and uniqueness of the open up mapping}
\label{sect:openup}

In this section, we prove Theorem~\ref{thm:openup}
using the theory of Riemann surfaces.
For background on Riemann surfaces and the
Riemann--Roch theorem, we refer to the books~\cite{Schlag}
or~\cite{Forster}.
If $\Gamma$ is a Jordan curve, denote the bounded and unbounded components of 
$\bC \setminus \Gamma$ by $\interior(\Gamma)$ and $\exterior(\Gamma)$,
respectively.

\begin{proof}[Proof of the existence in Theorem~\ref{thm:openup}]
Let $\gamma$ be a Jordan arc in $\C$.  Every point of $\gamma$ is accessible 
from $\C \setminus \gamma$; see~\cite[p.~164]{Newman1951}.
We will show that every point of $\gamma$ that is not an endpoint of $\gamma$ 
gives rise to two distinct accessible boundary points of $\C_\infty \setminus 
\gamma$, while an endpoint of $\gamma$ gives rise to one accessible 
boundary point of $\C_\infty \setminus \gamma$.
Recall that a boundary point $z$ of a domain $D$ is \emph{accessible}
from $D$ if there exists a Jordan arc $\ell$ with one endpoint at $z$ and
otherwise contained in $D$;
see, e.g.,~\cite[p.~162]{Newman1951} or~\cite[Ch.~II, \S 3, p.~35]{Goluzin}.
Such a Jordan arc $\ell$ is also called an \emph{end-cut}~\cite[p.~118]{Newman1951}.
Following Goluzin~\cite[pp.~36--37]{Goluzin}, two accessible boundary points 
$z_1$ and $z_2$ in $\partial D$ are regarded as distinct, if either $z_1 \neq 
z_2$ or if $z_1 = z_2$ but given two end-cuts $\ell_1, \ell_2$, there exists a 
neighborhood $U$ of $z_1$ such that $\ell_1, \ell_2$ cannot be joined in
$U \cap D$.
There exists a Jordan curve $\widetilde{\gamma}$ in $\C$ such that $\gamma$ is 
an arc of $\widetilde{\gamma}$, i.e., $\gamma \subseteq \widetilde{\gamma}$;
see~\cite[Ch.~VI, Thm.~14.5]{Newman1951} 
or~\cite[Cor.~17.23]{RaoStetkaerFournaisMoller2015}.
Let $\mathbb{D} = \{ z \in \C : \abs{z} < 1 \}$ be the unit disk.
By the Schoenflies theorem (see, e.g., \cite[Cor.~2.9]{Pommerenke1992}), there 
exists a homeomorphism $f : \C \to \C$ such that $f(\widetilde{\gamma}) = 
\partial \mathbb{D}$ is the unit circle, $f(\interior(\widetilde{\gamma})) = 
\interior(\partial \mathbb{D}) = \mathbb{D}$ and 
$f(\exterior(\widetilde{\gamma})) = \exterior(\partial \mathbb{D})$.
It is not hard to see that every point $z \in \partial \mathbb{D}$ gives rise 
to one accessible boundary point in $\mathbb{D}$ and to one in $\exterior(\partial \mathbb{D})$.
By the homeomorphism $f$, also every $z \in 
\widetilde{\gamma}$ gives rise to one accessible boundary point in 
$\interior(\widetilde{\gamma})$, denoted by $z^+$, and to one
accessible boundary point in $\exterior(\widetilde{\gamma})$, denoted by 
$z^-$.
In the domain $\C_\infty \setminus \gamma$, a point $z \in \gamma$ that 
is not an endpoint of $\gamma$ gives rise to two accessible 
boundary points:
$z^+$ accessible from $\interior(\widetilde{\gamma})$ and
$z^-$ accessible from $\exterior(\widetilde{\gamma})$.
An endpoint $z$ of $\gamma$ yields one accessible boundary point
($z^+ = z^-$ in $\C_\infty \setminus \gamma$), which follows 
from~\cite[Thm.~14.2, p.~162]{Newman1951}.
Set $\gamma^+ \coloneq \{ z^+ : z \in \gamma \}$ and
$\gamma^- \coloneq \{ z^- : z \in \gamma \}$.
As point sets, $\gamma^+ = \gamma^- = \gamma$, but as
sets of accessible boundary points, $\gamma^+$ and $\gamma^-$ are two Jordan 
arcs that are disjoint except at their endpoints.

\begin{figure}[t]
{\centering
\includegraphics[width=0.99\linewidth]{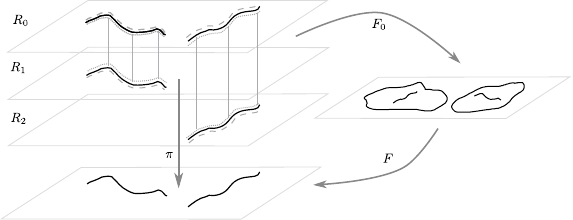}

}
\caption{The Riemann surface $R$ and the maps $F_0, \pi$, and $F$
in the proof of the existence in Theorem~\ref{thm:openup} 
in the case of two arcs.}
\label{riemannfig}
\end{figure}

Take $n+1$ copies of the Riemann sphere, denoted by $R_0, R_1, \ldots, R_n$. 
Cut $R_0$ along all the arcs $\gamma_1, \ldots, \gamma_n$.
For $j = 1, \ldots, n$, cut $R_j$ along $\gamma_j$ and join it crosswise to 
$R_0$ along the arc $\gamma_j$, i.e., for $z \in \gamma_j$ we identify
$z^+ \in \gamma_j^+$ in $R_0$ with $z^- \in \gamma_j^-$ in $R_j$, and
$z^- \in \gamma_j^-$ in $R_0$ with $z^+ \in \gamma_j^+$ in $R_j$.
(See, e.g., \cite[Ch.~4.3, p.~83]{Kalajdzievski2015} for the identification.)
This results in a Riemann surface which we denote by $R$; see, 
e.g.,~\cite[p.~8]{NasyrovBook} or the Russian translation of Hurwitz and
Courant~\cite[pp.~383--384 and p.~579]{HurwitzCourant1968}.
Figure~\ref{riemannfig} illustrates the construction of $R$ for $n=2$.
Let us sketch why $R$ is a Riemann surface.
The identification yields a topological surface,
on which we define the following charts:
$w = z$ in a neighborhood of a finite point that is not
a branch point (not an endpoint of an arc $\gamma_1, \ldots, \gamma_n$),
$w = \sqrt{z - \eta_k}$ in a neighborhood of a branch
point $\eta_k$ (an endpoint of one of the arcs $\gamma_1, \ldots, \gamma_n$),
and $w = 1/z$ in a neighborhood of $z = \infty$.

For each $j = 1, \ldots, n$, the above identification of $\gamma_j^+$ in $R_0$
with $\gamma_j^-$ in $R_j$ yields a simple arc $\hat{\gamma}_j^+$ in $R$,
and the identification of $\gamma_j^-$ in $R_0$ with $\gamma_j^+$ in $R_j$
yields a simple arc $\hat{\gamma}_j^-$ in $R$.
Since $\hat{\gamma}_j^+$ and $\hat{\gamma}_j^-$ are disjoint
except for their endpoints, $\hat{\gamma}_j \coloneq \hat{\gamma}_j^+ 
\cup \hat{\gamma}_j^-$ is a simple closed curve in $R$.

Note that $R$ is simply connected (i.e., has genus $0$) and compact.
One of the corollaries of the Riemann--Roch theorem
(see e.g.~\cite[pp.~130--131]{Forster}) says that
there is a biholomorphic mapping $F_0$ from $R$ 
onto the Riemann sphere $\bC_\infty$.
(Alternatively, since $R$ is compact and simply connected, the uniformization 
theorem implies the existence of $F_0$.)
We choose $F_0$ such that $\infty \in R_0$ is mapped to $\infty$
(if needed, this can be achieved by postcomposition with a M\"obius
transformation).

Let $\pi : R \to \bC_\infty$ be the canonical projection from $R$ onto the 
Riemann sphere, i.e., with $\pi(w^{(k)})=w$, where $w\in \bC$, and $w^{(k)} \in 
R_k$ is above $w$.
Then $F \coloneq \pi \circ F_0^{-1} : \bC_\infty \to \bC_\infty$
is a meromorphic function and hence $F$ is a rational function
(see e.g.~\cite[p.~11]{Forster} or \cite[Thm.~3.5.8]{Wegert2012}).
Because of the projection, $F$ is an $(n+1)$-to-$1$ mapping and 
thus $F$ has degree $n+1$.
The poles of $F$ are the images under $F_0$ of $\infty \in R_j$,
$j = 0, \ldots, n$.  In particular, $F$ has $n+1$ distinct poles
and these are simple.
Since $F_0^{-1}(\infty) = \infty \in R_0$, also $F(\infty) = \infty$ and
$F$ is of type $(n+1,n)$.

Consider the simply connected domains
$R_j \setminus \gamma_j = \C_\infty \setminus \gamma_j$,
$j = 1, \ldots, n$, in $R$ and the $n$-connected domain
$R_0 \setminus \cup_{j=1}^n \gamma_j
= \C_\infty \setminus \cup_{j=1}^n \gamma_j$ in $R$.
Since $F_0$ is biholomorphic,
\begin{equation*}
G_j \coloneq F_0(R_j \setminus \gamma_j) \subseteq \C_\infty, \quad j = 1, 
\ldots, n,
\end{equation*}
are simply connected domains and
\begin{equation*}
G_\infty \coloneq F_0(R_0 \setminus \cup_{j=1}^n \gamma_j) \subseteq \C_\infty
\end{equation*}
is an $n$-connected domain, and $\infty \in G_\infty$ by the above choice of $F_0$.
The domains $G_\infty, G_1, \ldots, G_n$ are disjoint since the sheets $R_0, 
R_1, \ldots, R_n$ are disjoint and $F_0$ is biholomorphic.
In particular, $G_1, \ldots, G_n$ are bounded.

By construction, the boundary of $R_j \setminus \gamma_j$ in $R$ is 
the simple closed curve $\hat{\gamma}_j$, and the boundary of
$R_0 \setminus \cup_{j=1}^n \gamma_j$ in $R$ consists of
the simple closed curves $\hat{\gamma}_1, \ldots, \hat{\gamma}_n$.
Since $F_0 : R \to \bC_\infty$ is biholomorphic,
$\Gamma_j \coloneq F_0(\hat{\gamma}_j)$, $j = 1, 
\ldots, n$, are disjoint Jordan curves in $\C_\infty$.
Moreover, $\Gamma_j$ is the boundary of $G_j$ and one boundary component of 
$G_\infty$.
Together, we obtain that $\partial G_\infty = \cup_{j=1}^n \Gamma_j$ consists 
of $n$ Jordan curves, and 
\begin{equation*}
K \coloneq \C_\infty \setminus G_\infty
= \cup_{j=1}^n (G_j \cup \Gamma_j) = \cup_{j=1}^n \overline{G}_j
\end{equation*}
is a compact set with $n$ components, each of which is the closure of a 
simply connected Jordan domain.
By construction, $F = \pi \circ F_0^{-1} : \C_\infty \setminus K \to \C_\infty \setminus \cup_{j=1}^n \gamma_j$ is bijective and conformal.  Similarly,
$F: G_j \to \bC_\infty \setminus \gamma_j$ is bijective and conformal.
This completes the proof.
\end{proof}

\begin{remark}
If the Jordan arcs $\gamma_j$ are $C^{k+}$ smooth, 
then the Jordan curves making up $\partial K$
are $C^{k+}$ smooth too, see~\cite[p.~206]{Widom},
where $C^{k+}$ means $k$ times continuously differentiable and
the $k$-th derivative is Lipschitz $\alpha$ for some $\alpha>0$.
Moreover, analyticity is also preserved,
that is, if $\gamma_j$ are analytic Jordan arcs,
then $\partial K$ consists of analytic Jordan curves,
see~\cite[p.~879]{KalmykovNagy2015}.
Both assertions are also clear from the proof of the existence in
Theorem~\ref{thm:openup}, since $\Gamma_j = F_0(\hat{\gamma}_j)$
is the image of a simple closed curve under a biholomorphic map.
\end{remark}

Next, we show that the  critical values of an open up mapping as in
Theorem~\ref{thm:openup} are precisely the endpoints of the arcs.

\begin{proposition} \label{prop:Prob1vsProb2} \label{prop:endpts_are_critvals}
Let $\gamma_1,\ldots,\gamma_n$ be disjoint 
Jordan arcs in the complex plane.
Denote the endpoints of $\gamma_j$ by $\cv_{2j-1}$ and $\cv_{2j}$,
$j=1,2,\ldots,n$.
If $F$ is a rational open up mapping of type $(n+1, n)$ for
$\gamma_1, \ldots, \gamma_n$, then the set of critical values of $F$ is
$\{ \cv_1, \ldots, \cv_{2n} \}$.
Moreover, $F^{-1}(\{ \cv_j \})$ consists of $n$ distinct points,
one of them has multiplicity two and is a critical point of $F$,
and the others have multiplicity one.
\end{proposition}

\begin{proof}
Let $F$ be a rational open up mapping of type $(n+1, n)$ for $\gamma_1,
\ldots, \gamma_n$.  In particular, $F : \bC_\infty \setminus K \to \bC_\infty
\setminus \cup_{j=1}^n \gamma_j$ is bijective and conformal.
Denote the boundary curves of $K$ by $\Gamma_1, \ldots, \Gamma_n$, labelled
such that $F(\Gamma_j) = \gamma_j$ for $j = 1, 2, \ldots, n$.

Fix $k \in \{ 1, 2, \ldots, 2n \}$ and let $\cv_k$ be an endpoint of $\gamma_j$
(i.e.~$k=2j-1$ or $k=2j$).
Since $F$ is of type $(n+1, n)$, there are $n+1$ pre-images of $\cv_k$ under
$F$ in $\bC$:
There is exactly one in $\interior(\Gamma_\ell)$ for each $\ell \neq j$,
and there are two on $\Gamma_j$.
(The pre-images cannot be in $\interior(\Gamma_j)$, since $F :
\interior(\Gamma_j) \to \bC_\infty \setminus \gamma_j$.)
Let $\cp_k \in \Gamma_j$ with $F(\cp_k) = \cv_k$.
Since $\cv_k$ are the endpoints of arcs, we must have $F'(\cp_k) = 0$
(otherwise $F$ is locally bijective at $\cp_k$ and $\cv_k = F(\cp_k)$ is an interior point of $\gamma_j$),
hence $\cv_k$ is a critical value of $F$ and $\cp_k$ is a double pre-image of
$\cv_k$ under $F$.
The set of critical values of $F$ is $\{ \cv_1, \ldots, \cv_{2n} \}$, since
$F$ is of type $(n+1, n)$ and cannot have any further critical values.
\end{proof}

Next, we show the uniqueness up to a linear transformation in
Theorem~\ref{thm:openup}, thus completing the proof.

\begin{proof}[Proof of the uniqueness in Theorem~\ref{thm:openup}]
Let $F, \widetilde{F}$ be rational functions that are open up mappings. 
We construct an analytic map $\varphi : \bC_\infty \to \bC_\infty$ with 
$\widetilde{F}(z) = F(\varphi(z))$ and show that $\varphi(z) = az + b$.

Let $G_\infty \coloneq \bC_\infty \setminus K$ be the region that
is mapped by $F$ onto $\bC_\infty \setminus \cup_{j=1}^n \gamma_j$.
For $j = 1, \ldots, n$, denote by $\Gamma_j$ the boundary curve of $K$ that is 
mapped by $F$ onto $\gamma_j$, and let $G_j = \interior(\Gamma_j)$.
Then $F : G_j \to \bC_\infty \setminus \gamma_j$ for $j = 1, \ldots, n$ 
and
$F : G_\infty \to \bC_\infty \setminus \cup_{j=1}^n \gamma_j$ are conformal 
(and bijective).
Introduce the same notation for $\widetilde{F}$, but with tildes.
Then $\varphi \coloneq F^{-1} \circ \widetilde{F} : \widetilde{G}_j \to G_j$
is conformal (and bijective) for all $j = 1, \ldots, n, \infty$, and
$\widetilde{F}(z) = F(\varphi(z))$ for $z \in \widetilde{G}_1 \cup \ldots \cup 
\widetilde{G}_n \cup \widetilde{G}_\infty$.

We extend $\varphi$ to an analytic function on $\bC$ with a simple pole at 
$\infty$.
Since $\widetilde{G}_1, \ldots, \widetilde{G}_n$ are Jordan regions, $\varphi$ 
extends to a homeomorphism $\widetilde{G}_j \cup \widetilde{\Gamma}_j \to 
G_j \cup \Gamma_j$ by the Osgood--Carath\'eodory theorem; 
see~\cite[Thm.~5.10e]{Henrici1974} or~\cite[Thm.~6.5.1]{Wegert2012}.
This is also true for $\widetilde{G}_\infty$ by subdividing the relevant parts 
in Jordan regions, similar to~\cite[p.~385]{Henrici1974}.

Let $j \in \{ 1, \ldots, n \}$ and $\widetilde{z}_0 \in \widetilde{\Gamma}_j$.
Then $\varphi : \widetilde{G}_j \cup \widetilde{\Gamma}_j \to G_j \cup 
\Gamma_j$ 
maps $\widetilde{z}_0$ to a point $z_0 \in \Gamma_j$.
First, consider $\widetilde{z}_0 \in \widetilde{\Gamma}_j$ that is not a 
critical point of $\widetilde{F}$, so that 
$\widetilde{F}(\widetilde{z}_0)$ is not an endpoint of $\gamma_j$
(see Proposition~\ref{prop:endpts_are_critvals}).
Then $z_0 \in \Gamma_j$ is not a critical point of $F$.
(Otherwise $F(z_0) = \widetilde{F}(\widetilde{z}_0)$ would be an endpoint of 
$\gamma_j$.)
Then there exist open neighborhoods $U$ of $z_0$, $\widetilde{U}$ of 
$\widetilde{z}_0$, and $V$ and $\widetilde{V}$ of $F(z_0)$, such that $F : U 
\to V$ and $\widetilde{F} : \widetilde{U} \to \widetilde{V}$ are conformal
(and bijective).
Without loss of generality, we have $V = \widetilde{V}$.
Then $F^{-1} \circ \widetilde{F} : \widetilde{U} \to U$ is analytic, 
maps $\widetilde{U} \cap \widetilde{G}_k \to U \cap G_k$ for $k \in \{ j, 
\infty \}$ and $\widetilde{U} \cap \widetilde{\Gamma}_j \to U \cap \Gamma_j$,
and coincides with $\varphi$ in $\widetilde{U} \cap \widetilde{G}_j$ and 
$\widetilde{U} \cap \widetilde{G}_\infty$.
Therefore, $\varphi$ extends to an analytic function in $\widetilde{U}$.

Next, let $\widetilde{z}_0 \in \widetilde{\Gamma}_j$ be a critical point of 
$\widetilde{F}$, then $z_0 \in \Gamma_j$ is a critical point of $F$.
By the above extension, $\varphi$ is analytic in a punctured neighborhood 
of $\widetilde{z}_0$ and continuous at $\widetilde{z}_0$, hence also analytic 
at $\widetilde{z}_0$.

Therefore, $\varphi$ is analytic in $\bC$ with a simple pole at $\infty$,
hence $\varphi(z) = a z + b$ with $a \neq 0$.
Then $\widetilde{F}(z) = F(\varphi(z)) = F(a z + b)$.
This shows uniqueness up to a linear transformation.
Finally, if $F$ and $\widetilde{F}$ both have the form $z + O(1/z)$
at infinity, then $\varphi(z) = z$ and $F = \widetilde{F}$.
\end{proof}

The open up mapping in Theorem~\ref{thm:openup}
depends only on the endpoints of the arcs and the topology of
$\bC \setminus \cup_{j=1}^n \gamma_j$, but does not depend on the specific
shape of the arcs.
This is shown in the next theorem, which is formulated for one arc,
but can be applied iteratively to allow deformation of all arcs.

\begin{theorem} \label{thm:deformation_of_arcs}
Let $\gamma_1, \ldots, \gamma_n$ be disjoint Jordan arcs in the
complex plane, and let $F : \bC_\infty \setminus K \to \bC_\infty \setminus
\cup_{j=1}^n \gamma_j$ be an open up mapping as in Theorem~\ref{thm:openup}.
Let $\widetilde{\gamma}_1$ be a Jordan arc with same endpoints as
$\gamma_1$ which is homotopic with fixed endpoints to $\gamma_1$ in
$\bC \setminus (\gamma_2 \cup \ldots \cup \gamma_n)$.  Then $F$ is also
an open up mapping for $\widetilde{\gamma}_1, \gamma_2, \ldots, 
\gamma_n$, i.e., there exists a compact set $\widetilde{K}$ bounded by $n$ 
disjoint Jordan curves such that
$F : \bC_\infty \setminus \widetilde{K} \to \bC_\infty \setminus
(\widetilde{\gamma}_1 \cup \gamma_2 \cup \ldots \cup \gamma_n)$ is conformal 
and bijective.
\end{theorem}

\begin{proof}
Since $F : \bC_\infty \setminus K \to \bC_\infty \setminus \cup_{j=1}^n 
\gamma_j$ is an open up mapping, $\partial K$ consists of $n$ disjoint Jordan 
curves $\Gamma_1, \ldots, \Gamma_n$, 
labeled such that $F(\Gamma_j) = \gamma_j$ for $j = 1, 2, \ldots, n$.
By Proposition~\ref{prop:endpts_are_critvals},
the endpoints $\cv_1, \cv_2$ of $\gamma_1$ are critical values of $F$.
Let $\cp_1, \cp_2 \in \Gamma_1$ be the critical points of $F$ with
$F(\cp_j) = \cv_j$ for $j = 1, 2$.

Let $f_0 : \cc{0, 1} \to \widetilde{\gamma}_1$, $f_1 : \cc{0, 1} \to \gamma_1$, 
be continuous and bijective functions
with $f_0(0) = f_1(0) = \cv_1$, $f_0(1) = f_1(1) = \cv_2$.
By assumption of the theorem, 
$f_0, f_1$ are homotopic with fixed endpoints in
$D = \bC \setminus (\gamma_2 \cup \ldots \cup \gamma_n)$, i.e., there exists a 
continuous function
$H : \cc{0, 1} \times \cc{0, 1} \to D$ with $H(0, t) = f_0(t)$ and $H(1, t) = 
f_1(t)$ for all $t \in \cc{0, 1}$ and $H(s, 0) = \cv_1$ and $H(s, 1) = 
\cv_2$ for all $s \in \cc{0, 1}$.
Then $H(\cc{0, 1} \times \cc{0, 1}) \subseteq D$ is compact.

Let $\gamma$ be a positively oriented Jordan curve in $D$ such that $H(\cc{0, 
1} \times \cc{0, 1}) \subseteq \interior(\gamma)$ and $\gamma_2, \ldots, 
\gamma_n \subseteq \exterior(\gamma)$.
Let $\Gamma^{(e)}$ be the pre-image in $\bC \setminus K$ of $\gamma$ under $F$.
Then $\Gamma^{(e)}$ is a positively oriented Jordan curve with $\Gamma_1$ in its 
interior (i.e., $\Gamma_1 \subseteq \interior(\Gamma^{(e)})$),
since $F : \bC_\infty \setminus K \to \bC_\infty \setminus \cup_{j=1}^n
\gamma_j$ is conformal and bijective.
Since $F : \interior(\Gamma_1) \to \bC_\infty \setminus \gamma_1$ is also 
conformal and bijective, there exists a Jordan curve $\Gamma^{(i)} \subseteq 
\interior(\Gamma_1)$ such that $F : \Gamma^{(i)} \to \gamma$ is bijective.
Note that $\Gamma^{(i)}$ is negatively oriented.

Let $A$ denote the (open) ring-domain bounded by $\Gamma^{(e)}$ and $\Gamma^{(i)}$.
Then $F$ is holomorphic on $A$, and $\Gamma_1 \subseteq A$.
For $w \in \interior(\gamma)$, the winding of $F-w$
along $\Gamma^{(i)} \cup 
\Gamma^{(e)}$ is $W(F - w; \Gamma^{(i)} \cup \Gamma^{(e)}) = 2$, hence
$F : A \to \interior(\gamma)$ is $2$ to $1$ by the argument principle.

We have $\widetilde{\gamma}_1 \subseteq \interior(\gamma)$ by the definition of 
$\gamma$.
Hence the pre-image of $\widetilde{\gamma}_1$ under $F$ in $A$ consists of two 
Jordan arcs connecting $\cp_1, \cp_2$.
Since $F$ is conformal in $A \setminus \{ \cp_1, \cp_2 \}$, the two arcs 
cannot intersect except at $\cp_1, \cp_2$, and hence form a Jordan curve 
$\widetilde{\Gamma}_1$.
We orient $\widetilde{\Gamma}_1$ in the negative sense.
We then have for $w \in \interior(\gamma) \setminus \widetilde{\gamma}_1$ that 
$W(F-w; \Gamma^{(e)} \cup \widetilde{\Gamma}_1) = 1$, hence $F$ maps the ring 
domain bounded by $\widetilde{\Gamma}_1$ and $\Gamma^{(e)}$ bijectively onto 
$\interior(\gamma) \setminus \widetilde{\gamma}_1$.
This implies that
$F : \bC_\infty \setminus \widetilde{K} \to \bC_\infty \setminus ( 
\widetilde{\gamma}_1 \cup \cup_{j=2}^n \gamma_j )$
is conformal and bijective and hence an open up mapping, where $\bC_\infty 
\setminus \widetilde{K}$ is the unbounded domain with boundary 
$\widetilde{\Gamma}_1 \cup \cup_{j=2}^n \Gamma_j$.
\end{proof}

We saw in Proposition~\ref{prop:endpts_are_critvals} that an open up mapping
has the endpoints of the arcs as critical values.
Theorem~\ref{thm:deformation_of_arcs} clarifies the difference between
an open up mapping and rational functions with critical values
at the endpoint of the arcs.
While the information on the critical values is present in both problems,
the difference is the additional
``topological'' 
information about the arcs
in Theorem~\ref{thm:openup}, which is not present in 
Theorem~\ref{thm:critvals}.
We give an example for this difference
in Section~\ref{sect:example}.

\section{Rational functions with prescribed critical values}
\label{sect:critvals}

In this section, we first prove Theorem~\ref{thm:critvals}.
Afterwards, we also consider an equivalent polynomial formulation, which can be 
suitable for the computation of rational functions with prescribed critical 
values.

The existence of solutions in Theorem~\ref{thm:critvals} readily follows from 
Theorem~\ref{thm:openup} and Proposition~\ref{prop:endpts_are_critvals}.
The exact number of solutions is derived with results on ramified coverings of 
the Riemann sphere.

\begin{proof}[Proof of Theorem~\ref{thm:critvals}.]
\emph{Step 1: Existence.}
Connect the points $\cv_1, \ldots, \cv_{2n}$ pairwise by Jordan arcs
that do not intersect each other.  By Theorem~\ref{thm:openup} there exists
a rational open up mapping $F$ of type $(n+1, n)$, which by 
Proposition~\ref{prop:Prob1vsProb2} is a
solution in Theorem~\ref{thm:critvals}.

\emph{Step 2: Number of normalized solutions.}
Let $F$, $F_1$ be two rational functions of type $(n+1, n)$ with critical
values $\cv_1, \ldots, \cv_{2n}$ (we do not impose the normalization at infinity yet).
Following Mednykh~\cite{Mednykh}, we call $F, F_1$ equivalent, if there exists
a homeomorphism $\varphi : \C_\infty \to \C_\infty$ with $F = F_1 \circ 
\varphi$.
It follows that $\varphi$ is a M\"obius transformation.

By assumption, the critical values $\cv_1, \ldots, \cv_{2n} \in \C$ are
distinct and hence are ramification points of order $2$ of the rational
functions.  In this case, the number of non-equivalent coverings
(i.e., of equivalence classes) in the sense of Mednykh is given by the Hurwitz 
numbers $H_n$ in~\eqref{eqn:hurwitz_numbers},
see~\cite{Mednykh} and Hurwitz' original article~\cite[p.~22]{Hurwitz1891};
see also~\cite[p.~290]{LandoZvonkin} and~\cite[Eqn.~(4.10)]{CrescimannoTaylor}.

In the second step, we estimate the number of normalized functions in
Theorem~\ref{thm:critvals} in each equivalence class.
Note first that, by Remark~\ref{rem:lineartransformation}, each function
is equivalent to one that is normalized at infinity.
Next, let $F, F_1$ be \emph{two equivalent} rational functions
of type $(n+1, n)$ with critical values $\cv_1, \ldots, \cv_{2n}$
and normalized at infinity by $z + O(1/z)$.  Let us write
\begin{equation*}
F(z) = z + \sum_{j=1}^n \frac{r_j}{z - p_j},
\end{equation*}
where $p_1, \ldots, p_n \in \C$ are distinct,
and $r_j \neq 0$ for $j = 1, \ldots, n$.
As shown above, $F = F_1 \circ \varphi$ with a M\"obius transformation
$\varphi$.  In particular, $p$ is a pole of $F$ if and only if $\varphi(p)$ is a 
pole of $F_1$.
In particular, $\varphi^{-1}(\infty)$ is a pole of $F$ and we distinguish two 
cases.

If $\varphi^{-1}(\infty) = \infty$, then $\varphi$ is a linear transformation
of the form $\varphi(z) = az + b$ with nonzero $a \in \C$.  The normalization 
of $F$ and $F_1$ at infinity and $F = F_1 \circ \varphi$ imply $\varphi(z) = 
z$, hence $F = F_1$.

Otherwise, $\varphi^{-1}(\infty)$ is a finite pole of $F$ and there exists $j_0 
\in \{ 1, \ldots, n \}$ such that $\varphi^{-1}(\infty) = p_{j_0} \in \C$.
Then $\varphi$ has the form
\begin{equation} \label{eqn:moebius}
\varphi(z) = \frac{a z + b}{z - p_{j_0}}, \quad
\varphi^{-1}(z) = \frac{p_{j_0} z + b}{z - a},
\end{equation}
with $a, b \in \C$ and $- a p_{j_0} - b \neq 0$.  We compute
\begin{align}
F_1(z) = F(\varphi^{-1}(z))
&= \frac{p_{j_0} z + b}{z - a} + \sum_{j=1}^n r_j \frac{z - a}{(p_{j_0} z + b) 
- p_j (z-a)} \notag \\
&= p_{j_0} + \frac{a p_{j_0} + b}{z - a} + \sum_{j=1}^n r_j \frac{z - 
a}{(p_{j_0} - p_j) z + a p_j + b}. \label{eqn:F1_in_finiteness}
\end{align}
We distinguish the cases $j = j_0$ and $j \neq j_0$ in the sum.
If $j \neq j_0$ then
\begin{equation*}
\frac{z - a}{(p_{j_0} - p_j) z + a p_j + b}
= \frac{1}{p_{j_0} - p_j} \cdot \frac{z - a}{z - \varphi(p_j)}
= \frac{1}{p_{j_0} - p_j} + \frac{\frac{\varphi(p_j) - a}{p_{j_0} - p_j}}{z - 
\varphi(p_j)}.
\end{equation*}
Inserting this in~\eqref{eqn:F1_in_finiteness} yields
\begin{equation*}
F_1(z)
= \frac{r_{j_0}}{a p_{j_0} + b} (z - a) + p_{j_0} 
+ \frac{a p_{j_0} + b}{z - a}
+ \sum_{j \neq j_0} \frac{r_j}{p_{j_0} - p_j}
+ \sum_{j \neq j_0} \frac{r_j \frac{\varphi(p_j) - a}{p_{j_0} - p_j}}{z - 
\varphi(p_j)}.
\end{equation*}
By assumption, $F(z) = z + O(1/z)$ at $\infty$.  Comparing the coefficient of 
$z$ yields
\begin{equation} \label{eqn:normalize_at_infty}
\frac{r_{j_0}}{a p_{j_0} + b} = 1,
\end{equation}
and comparing the constant coefficient yields
\begin{equation} \label{eqn:normalize_at_infty_cst}
a = p_{j_0} + \sum_{j \neq j_0} \frac{r_j}{p_{j_0} - p_j}.
\end{equation}
Moreover, by~\eqref{eqn:normalize_at_infty}, $b = r_{j_0} - a p_{j_0}$.
Thus, $a$ and $b$ are uniquely determined by $F$ and $p_{j_0}$.
This shows that $\varphi$ is fully determined by $F$ and a choice
of one of the $n+1$ poles of $F$.
Thus, each of the $H_n$ many equivalence classes contains at most $n+1$ 
distinct normalized functions, which establishes the \emph{bound} $(n+1) H_n$ on 
the number of normalized rational functions in Theorem~\ref{thm:critvals}.
We proceed to derive the \emph{exact} number.

If $n = 1$, there is one equivalence class since $H_1 = 1$.  Thus the bound 
gives $2$, but there is only one solution in Theorem~\ref{thm:critvals}, 
which can be seen as follows;
see also Proposition~\ref{prop:casen1} where the solution is derived
explicitly.  If $F(z) = z + \frac{r_1}{z - p_1}$ and if we choose
$j_0 = 1$, then $a = p_1$ by~\eqref{eqn:normalize_at_infty_cst} and thus
$F_1(z) = z + \frac{r_1}{z - a} = F(z)$, so that there is only
a single normalized solution in the only equivalence class.

If $n \geq 2$, the $n+1$ normalized functions in each equivalence class are 
distinct, which we show next.
Consider as above $F = F_1 \circ \varphi$ with a M\"obius transformation 
$\varphi$ of the form~\eqref{eqn:moebius} for some $j_0 \in \{ 1, \ldots, n 
\}$.  We show by contradiction that $F$ and $F_1$ are distinct, and therefore 
assume that $F = F_1$.  This leads to
\begin{equation} \label{eqn:F_selftransform}
F(z) = (F \circ \varphi)(z).
\end{equation}
In particular, $\varphi(\{ p_1, \ldots, p_n, p_\infty = \infty \}) = \{ p_1, 
\ldots, p_n, p_\infty = \infty \}$, that is $\varphi$ permutes the poles of 
$F$.
Since the set of poles of $F$ is finite, there exists a positive integer $N$ 
such that $\varphi^N$, the $N$-th iteration of $\varphi$, satisfies 
$\varphi^N(p_j) = p_j$ for $j = 1, \ldots, n, \infty$, hence the M\"obius 
transformation $\varphi^N$ has $n+1 \geq 3$ fixed points and thus is the 
identity, $\varphi^N = \id$.
Note that $\varphi \neq \id$, since $\varphi$ has a finite pole at $p_{j_0} 
\in \C$.

Next, we show that $\varphi$ has two distinct fixed points.
Since $\varphi \neq \id$, it has one or two fixed points.
Let us assume that $\varphi$ has only one fixed point $A \in \C_\infty$, we 
shall reach a contradiction.
Let $z = \chi(u)$ be a M\"obius transformation with $\chi(\infty) = A$.
The M\"obius transformation $\chi^{-1} \circ \varphi \circ \chi$ fixes infinity 
and has only one fixed point, hence 
$(\chi^{-1} \circ \varphi \circ \chi)(u) = u + \beta$ for some $\beta \in \C$.
Moreover, the set $\chi^{-1}(\{ p_1, \ldots, p_n, \infty \})$ is invariant 
under $\chi^{-1} \circ \varphi \circ \chi$, 
hence $\beta = 0$ and $\varphi = \id$, a contradiction.
Thus, $\varphi$ has two distinct fixed points $A, B \in \C_\infty$, $A \neq B$.

There exists a M\"obius transformation $z = \psi(w)$ such that $\psi(0) = A$ 
and $\psi(\infty) = B$.
Therefore $\psi^{-1}( \varphi(\psi(w))) = \lambda w $ for some $\lambda \in 
\bC$, where $\lambda \neq 1$ since $\varphi \neq \id$.
Transforming $\varphi$ to the $w$-plane, we can write 
\begin{equation*}
w = \psi^{-1} \bigl( \varphi^N(\psi(w)) \bigr)
= (\psi^{-1} \circ \varphi \circ \psi)^N(w)
= \lambda^N w
\end{equation*}
which implies that $\lambda^N = 1$.  Since $\lambda \neq 1$, we have $N \geq 
2$, and we can choose the smallest number $\nu \geq 2$ with $\lambda^\nu = 1$.
We transform $F$ to the $w$-plane and define $G(w) \coloneq F(\psi(w))$.
Then~\eqref{eqn:F_selftransform} yields $G(w) = F(\psi(w)) = 
F(\varphi(\psi(w))) = F(\psi(\lambda w)) = G(\lambda w)$.
Note that $F$ and $G$ have the same critical values $\cv_1, \ldots, \cv_{2n}$.  
Let $w_1, \ldots, w_{2n} \in \C$ be critical points of $G$ with $G(w_j) = 
\cv_j$ for $j = 1, \ldots, 2n$.
Since $\deg(G) = \deg(F) = n+1$, we have $\deg(G') \leq 2n+1$, i.e., $G$ has at 
most $2n+1$ critical points.
Now, $G'(w) = \lambda G'(\lambda w)$, implies that with $w_j$, also $\lambda 
w_j, \lambda^2 w_j, \ldots, \lambda^{\nu - 1} w_j$ are critical points 
corresponding to the critical value $\cv_j$.
If $w_j \neq 0$, these are pairwise distinct since $\nu \geq 2$ is minimal with 
the property $\lambda^\nu = 1$.
Since critical points corresponding to different critical values are distinct, 
we obtain that $G$ has strictly more than $2n+1$ critical points, a 
contradiction.
Therefore, the assumption $F = F_1$ is wrong, and there are indeed exactly 
$n+1$ distinct normalized rational functions in each equivalence class.  This 
completes the proof of Theorem~\ref{thm:critvals}.
\end{proof}

In view of computing the rational functions with prescribed critical values,
we reformulate Theorem~\ref{thm:critvals} in terms of polynomials and 
show the equivalence of the two formulations 
in Proposition~\ref{prop:Prob5vsProb2}.

\begin{theorem} \label{thm:critvals_poly}
Let $\cv_1, \cv_2 \ldots, \cv_{2n} \in \bC$ be distinct.
Then there exist polynomials
\begin{equation}
P(z) = z^{n+1} + \sum_{j=0}^{n-1} p_j z^j \quad \text{and} \quad
Q(z) = z^n + \sum_{j=0}^{n-1} q_j z^j,
\label{monic_PQ}
\end{equation}
i.e., with $p_{n+1} = q_n = 1$ and $p_n = 0$, and points
$\cp_1, \cp_2, \ldots, \cp_{2n} \in \bC$ such that
\begin{align}
P(\cp_j) - \cv_j Q(\cp_j) &= 0, \quad j = 1, 2, \ldots, 2n,
\label{poly_critvals} \\
P'(\cp_j) Q(\cp_j) - P(\cp_j) Q'(\cp_j) &= 0. \quad j = 1, 2, \ldots, 2n,
\label{poly_critpts}
\end{align}
and
\begin{equation}
Q(\cp_j) \neq 0, \quad j = 1, 2, \ldots, 2n. \label{poly_nopole}
\end{equation}
If $n = 1$, the solution $(P, Q)$ is unique.
If $n \geq 2$, the number of solutions $(P, Q)$ is
$(n+1) H_n$ with the Hurwitz number $H_n$ in~\eqref{eqn:hurwitz_numbers}.
\end{theorem}

There are $4n$ equations in~\eqref{poly_critvals}
and~\eqref{poly_critpts} for the
$4n$ unknowns $p_{n-1},\ldots,p_1,p_0$, $q_{n-1},\ldots,q_1$, $q_0$
and $\cp_1, \cp_2, \ldots, \cp_{2n}$ in Theorem~\ref{thm:critvals_poly},
which is twice as many as in Theorem~\ref{thm:critpts_poly}.

Equation~\eqref{poly_critvals} prescribes the values of $F = P/Q$ at the
points $\cp_1, \ldots, \cp_{2n}$, provided that~\eqref{poly_nopole} holds,
and $\cp_1, \ldots, \cp_{2n}$ are the critical points of $F$ 
by~\eqref{poly_critpts}.
If $P, Q$ in~\eqref{monic_PQ} satisfy~\eqref{poly_critvals}
and~\eqref{poly_critpts} but
$Q(\cp_j) = 0$ for some $j$, then also $P(\cp_j) = 0$ and $F = P/Q$
does not satisfy Theorem~\ref{thm:critvals}
since the degree of $F$ is too small; see
Proposition~\ref{prop:minimal}.

We prove Theorem~\ref{thm:critvals_poly} by showing that the solutions $(P, Q)$ 
correspond to the rational functions $F = P/Q$ in Theorem~\ref{thm:critvals}.
The different normalization of $F$ and $(P, Q)$ is not essential; see 
Remark~\ref{rem:lineartransformation}.

\begin{proposition} \label{prop:Prob5vsProb2}
Let $\cv_1, \ldots, \cv_{2n} \in \bC$ be distinct.
\begin{enumerate}
\item \label{item:Prob5vsProb2_1}
If $(P, Q)$ satisfy the equations in Theorem~\ref{thm:critvals_poly},
then $F = P/Q$
is a rational function of type $(n+1, n)$ with critical values
$\cv_1, \ldots, \cv_{2n}$, i.e., $F$ satisfies Theorem~\ref{thm:critvals},
$P$ and $Q$ are coprime
and the points $\cp_1, \ldots, \cp_{2n}$ are distinct.

\item \label{item:Prob5vsProb2_2}
If $F$ is a function in Theorem~\ref{thm:critvals},
then for all $a, b \in \bC$ with $a \neq 0$, also $F(az+b)$ is a 
function as in Theorem~\ref{thm:critvals}
and there exist unique $a, b$ such that
$F(az+b) = z + O(1/z)$ for $z \to \infty$.  Moreover, 
there exist $a, b$ such that $F(az+b) = P(z)/Q(z)$ and $(P, Q)$ 
satisfy Theorem~\ref{thm:critvals_poly}.
\end{enumerate}
\end{proposition}

\begin{proof}
Part~\ref{item:Prob5vsProb2_1} is obvious.  For part~\ref{item:Prob5vsProb2_2}, 
let $F$ be
of type $(n+1, n)$ such that
\begin{equation*}
F(\cp_j) = \cv_j, \quad F'(\cp_j) = 0, \quad j = 1, 2, \ldots, 2n.
\end{equation*}
Let $G(z) = F(a z + b)$, with $a, b \in \bC$ and $a \neq 0$, then
\begin{equation*}
G((\cp_j-b)/a) = \cv_j, \quad G'((\cp_j-b)/a) = F'(\cp_j) a = 0, \quad j = 1, 2, 
\ldots, 2n,
\end{equation*}
and $G$ also is as in Theorem~\ref{thm:critvals}.

If $F = P/Q$, then $(P, Q)$
satisfies~\eqref{poly_critvals}--\eqref{poly_critpts}
and $\deg(P) = n+1$ and $\deg(Q) = n$ (see Proposition~\ref{prop:minimal}),
but the normalization \eqref{monic_PQ} is not necessarily satisfied.
However, by a suitable affine transformation of the argument,
$G(z) = F(az+b)$, we obtain $G = P_1/Q_1$, such that $(P_1, Q_1)$
satisfies~\eqref{monic_PQ}.
Moreover, $Q_1(\cp_j) \neq 0$ (otherwise, $P_1(\cp_j) = 0$
by~\eqref{poly_critvals}, and $G$ would not be of type $(n+1, n)$).
Hence $(P_1, Q_1)$ is a solution in Theorem~\ref{thm:critvals_poly}.
The remaining assertions are clear.
\end{proof}

\section{Rational functions with prescribed critical points}
\label{sect:critpts}

We first reformulate the task of finding a rational function
$F = P/Q$ with prescribed critical points in terms of the polynomials $P$ and 
$Q$, and show the equivalence of the two formulations.

\begin{theorem} \label{thm:critpts_poly}
Let $\cp_1,\cp_2,\ldots,\cp_{2n} \in \bC$ be distinct.
Then there exist polynomials
\begin{equation}
P(z) = \sum_{j=0}^{n+1} p_j z^j
\quad \text{and} \quad
Q(z) = \sum_{j=0}^n q_j z^j,
\label{PQform}
\end{equation}
where $p_0, p_1, \ldots, p_{n+1}, q_0, q_1, \ldots, q_n \in \bC$
and $p_{n+1} \neq 0$ and $q_n \neq 0$, such that
\begin{equation}
\label{homPQs}
P'(\cp_j)Q(\cp_j)-P(\cp_j)Q'(\cp_j) =0,
\quad j=1,\ldots,2n.
\end{equation}
Moreover, each solution $(P, Q)$ can be normalized by
\begin{equation}
p_{n+1} = q_n = 1 \quad \text{and} \quad p_n = 0.
\label{leadingcoeffassump}
\end{equation}
The number of normalized solutions is bounded from above by the
Catalan number $C_n = \frac{1}{n+1} \binom{2n}{n}$.
\end{theorem}

\begin{remark} \label{rem:critpts_poly}
Each solution in Theorem~\ref{thm:critpts_poly} yields infinitely many
solutions, which can be seen as follows.
Let $(P, Q)$ be a solution in Theorem~\ref{thm:critpts_poly}, i.e., $P$ and $Q$
are polynomials of degrees $\deg(P) = n+1$ and $\deg(Q) = n$ that
satisfy~\eqref{homPQs}.  Then, for each $c, d \in \bC \setminus \{ 0 \}$
and $p \in \bC$, also $(c P + p Q, d Q)$ is a solution in
Theorem~\ref{thm:critpts_poly}.
We introduce an equivalence relation on the solutions: $(R, S) \sim (P, Q)$ if
and only if $R = cP+pQ$, $S = dQ$ for some $c, d \in \C \setminus \{ 0 \}$ and 
$p \in \C$.
In particular, the equivalence class of $(P, Q)$ contains a unique
representative satisfying~\eqref{leadingcoeffassump}, which is obtained
by $c = 1/p_{n+1}$, $d = 1/q_n$ and $p = - p_n/(q_n p_{n+1})$.
\end{remark}

First, we show that a rational function $F = P/Q$ in Theorem~\ref{thm:critpts}
corresponds to a solution $(P, Q)$ in Theorem~\ref{thm:critpts_poly}.

\begin{proposition} \label{prop:critpts_rat_vs_poly}
Let $\cp_1, \ldots, \cp_{2n} \in \bC$ be distinct.
Then $(P, Q)$ is a solution in Theorem~\ref{thm:critpts_poly} if, and only if,
$F = P/Q$ is a rational function as in Theorem~\ref{thm:critpts},
i.e., $F$ is of type $(n+1, n)$ with critical points $\cp_1, \ldots, \cp_{2n}$.
\end{proposition}

\begin{proof}
If $(P, Q)$ is a solution in Theorem~\ref{thm:critpts_poly},
then the polynomials $P$ and $Q$ are coprime, which can be seen as follows.
Assume to the contrary that $P, Q$ have a non-constant
common factor $C$ and write $P = C R$, $Q = C S$.
Then $P' Q - P Q' = C^2 (R' S - R S')$, which has $2n$ 
distinct zeros by~\eqref{homPQs}.
Since $C$ has at most $\deg(C) \geq 1$ distinct zeros,
$R'S-RS'$ has at least $2n-\deg(C)$ many distinct zeros.
On the other hand, $\deg(R'S-RS') = 2n - 2 \deg(C)$.
This implies $R'S-RS' \equiv 0$, and thus $R = c S$ with $c \in \bC \setminus \{ 0 \}$,
in contradiction to $\deg(P) = n+1 > \deg(Q) = n$.
Then $F = P/Q$ is of type $(n+1, n)$ and is a solution in
Theorem~\ref{thm:critpts}.
The converse is clear.
\end{proof}

Thus, Theorem~\ref{thm:critpts} and Theorem~\ref{thm:critpts_poly} are
equivalent.
In particular, the number of equivalence classes of solutions 
is the same in both theorems.
We show the existence of solutions for Theorem~\ref{thm:critpts_poly} using 
algebraic tools, in particular Hilbert's Nullstellensatz.
The bound on the number of solutions will be established for 
Theorem~\ref{thm:critpts} using a result of Goldberg.

\begin{proof}[Proof of Theorems~\ref{thm:critpts} and~\ref{thm:critpts_poly}.]
\emph{Step 1: Existence.}
We show that there are polynomials $P$ and $Q$ as in~\eqref{PQform} 
satisfying
\begin{equation}
P'(z) Q(z) - P(z) Q'(z) = \prod_{j=1}^{2n} (z - \cp_j). \label{eqn:PQ_aim}
\end{equation}
Then $\deg(P) = n+1$, $\deg(Q) = n$ and $(P, Q)$ is a solution in Theorem~\ref{thm:critpts_poly}.

We write
\begin{equation} \label{eqn:def_c_ell}
\prod_{j=1}^{2n} (z - \cp_j) = \sum_{\ell =0}^{2n} c_\ell z^\ell
\end{equation}
with constants $c_0, \ldots, c_{2n} \in \bC$ given by
$\cp_1, \ldots, \cp_{2n}$.  Note that $c_{2n} = 1$.
We also write
\begin{equation} \label{eqn:def_rho_ell}
P'(z) Q(z) - P(z) Q'(z) = \sum_{\ell = 0}^{2n} \rho_\ell z^\ell,
\end{equation}
where the coefficients are
\begin{align}
\rho_\ell &\coloneq
\sum_{\substack{j+k=\ell \\ j\le n, k\le n}} (j+1) p_{j+1} q_k
- \sum_{\substack{j+k=\ell \\ j\le n+1, k\le n-1}} (k+1) p_{j} q_{k+1}
\label{eqn:rho} \\
&\;= \sum_{\substack{j+k = \ell+1 \\ 0 \leq j \leq n+1, 0 \leq k \leq n}}
(j-k) p_j q_k.
\notag
\end{align}
In particular, $\rho_{2n} = p_{n+1} q_n$.
The structure of the coefficients $\rho_\ell$ is important and we will use the 
following fact later: For each $\ell$, 
$\rho_\ell$ is a 
homogeneous polynomial of order $2$ in the variables $p_0, p_1, \ldots, 
p_{n+1}$, $q_0, q_1, \ldots, q_n$ and contains only products $p_j q_k$ with
$j + k - 1 = \ell$.  Conversely, each product $p_j q_k$
(with $j\ne k$) appears only in 
$\rho_{j + k - 1}$ and with coefficient $j-k$.

Then~\eqref{eqn:PQ_aim} is equivalent to the system of polynomial equations
\begin{equation}
\rho_\ell = c_\ell, \quad \ell =0,1,\ldots,2n.
\label{overlinerhoeqs}
\end{equation}
To show that~\eqref{overlinerhoeqs} has a solution, 
we consider the ideal generated by $\rho_0 -c_0 ,\ldots$, $\rho_{2n}-c_{2n}$:
\begin{align*}
I &\coloneq \left\langle \rho_0 -c_0 ,\ldots,\rho_{2n}-c_{2n} \right\rangle 
\\
&= \bigg\{ \sum_{\ell = 0}^{2n} A_\ell  (\rho_\ell -c_\ell) :
A_\ell \in\bC[p_{n+1},\ldots,p_1,p_0,q_{n},\ldots,q_1,q_0],
\ell = 0,\ldots,2n \bigg\}.
\end{align*}
By the weak form of Hilbert's Nullstellensatz, 
\eqref{overlinerhoeqs} has a solution 
if and only if
\begin{equation*}
I \neq \bC[p_{n+1},\ldots,p_1,p_0,q_{n},\ldots,q_1,q_0];
\end{equation*}
see, e.g.,~\cite[Thm.~2.2.3]{AdamsLoustaunau}
or~\cite[Thm.~1, p.~177]{CLOIdeals}.
Hence, we are going to show that there exists a polynomial $Y \notin I$, i.e., 
such that
\begin{equation} \label{eqn:Y_notin_ideal}
Y = \sum_{\ell = 0}^{2n} A_\ell (\rho_\ell - c_\ell)
\end{equation}
does not hold for any $A_0, \ldots, A_{2n} \in
\bC[p_{n+1},\ldots,p_1,p_0,q_{n},\ldots,q_1,q_0]$.

Consider the equations in~\eqref{overlinerhoeqs} with nonzero right hand side
\[
L_1 \coloneq \{ \ell \in \{0,1,\ldots,2n\} : c_\ell \ne 0\}.
\]
Note that $\abs{L_1} > 1$.  Indeed, $c_{2n} = 1$ so $2n \in L_1$, and 
$\abs{L_1} = 1$ would imply $\prod_{j=1}^{2n} (z - \cp_j) = z^{2n}$, which 
contradicts that $\cp_1, \ldots, \cp_{2n}$ are distinct.
We determine the minimal element in these equations with respect to
the ordering
\begin{equation*}
p_{n+1} > p_n > p_{n-1}>\ldots>p_1>p_0>q_0>q_1>\ldots>q_{n-1} > q_n
\end{equation*}
of the variables with some order, e.g.,
the degree lexicographic order (also called graded lexicographic order, see, e.g.~\cite[p.~58]{CLOIdeals})
\begin{equation*}
p_{j_1} q_{k_1} = \min \{ p_r q_s : \coeff(p_r q_s, \rho_\ell) \neq 0 \text{ 
and } \ell \in L_1 \},
\end{equation*}
where $\coeff(p_r q_s, \rho_\ell)$ denotes the coefficient of $p_r q_s$ in 
$\rho_\ell$.
Let $\ell_1 = j_1 + k_1 - 1$ be the unique index such that $p_{j_1} q_{k_1}$ 
appears in $\rho_{\ell_1}$.  In particular, $\ell_1 \in L_1$.
Note that $j_1 \neq k_1$, otherwise the 
coefficient of $p_{j_1} q_{k_1}$ would be zero.

For $\ell \in L_1$, we write
\begin{equation*}
\rho_\ell - c_\ell = \rho_\ell - \frac{c_\ell}{c_{\ell_1}} \rho_{\ell_1}
+ \frac{c_\ell}{c_{\ell_1}} (\rho_{\ell_1} - c_{\ell_1}), \quad \ell \in L_1.
\end{equation*}
Then~\eqref{eqn:Y_notin_ideal} can be written as
\begin{equation}
\left.
\begin{aligned}
Y
&= \sum_{\ell \notin L_1} A_\ell \rho_\ell
+ \sum_{\ell \in L_1} A_\ell \frac{c_\ell}{c_{\ell_1}} (\rho_{\ell_1} - c_{\ell_1})
+ \sum_{\ell \in L_1} A_\ell \left( \rho_\ell - \frac{c_\ell}{c_{\ell_1}} \rho_{\ell_1} \right)
\\
&= \sum_{\ell \notin L_1} A_\ell \rho_\ell
+ \widehat{A}_{\ell_1} (\rho_{\ell_1} -c_{\ell_1})
+ \sum_{\ell \in L_1 \setminus \{ \ell_1 \}} A_\ell \left( \rho_\ell - \frac{c_\ell}{c_{\ell_1}} \rho_{\ell_1} \right),
\end{aligned}
\right. 
\label{elimdeq}
\end{equation}
where $\widehat{A}_{\ell_1} \coloneq \sum_{\ell \in L_1} A_\ell \frac{c_\ell}{c_{\ell_1}}$.
We will construct a polynomial 
$Y$ for which this equation has no solutions.

If $p_{j_1}$ appears in some of the terms in $\rho_\ell$, then it
appears only in the term $p_{j_1} q_{k(\ell)}$,
where $k(\ell) \coloneq \ell + 1 - j_1$.  
We define
\begin{align*}
L_2 &\coloneq \{ \ell \in \{ 0, 1, \ldots, 2n \} : p_{j_1} \text{ appears in } 
\rho_\ell \}
\\
&= \{ \ell \in \{ 0, 1, \ldots, 2n \} : \coeff(p_{j_1} q_{k(\ell)}, \rho_\ell) 
\neq 0 \}.
\end{align*}
Substituting $0$ for each variable in $\{ p_0, p_1, \ldots, 
p_{n+1} \} \setminus \{ p_{j_1} \}$,
we obtain (after this substitution) $\rho_\ell = 0$ for $\ell \notin L_2$,
and $\rho_\ell = (j_1 - k(\ell)) p_{j_1} q_{k(\ell)}$ for $\ell \in L_2$.
Denoting by $B_\ell$ the polynomial obtained from $A_\ell$
through the substitution,
\eqref{elimdeq} becomes
\begin{equation}
\begin{split}
Y= &\sum_{\ell \in L_2 \setminus L_1} B_\ell \, (j_1 - k(\ell)) p_{j_1} q_{k(\ell)}
+ B_{\ell_1} ( (j_1-k_1) p_{j_1} q_{k_1} - c_{\ell_1}) \\
&+ \sum_{\ell \in L_2 \cap L_1 \setminus \{\ell_1\}}
B_\ell \left( (j_1-k(\ell)) p_{j_1} q_{k(\ell)} - 
\frac{c_\ell}{c_{\ell_1}} (j_1-k_1) p_{j_1} q_{k_1} \right).
\end{split} \label{subselimdeq}
\end{equation}

If $L_2 \cap L_1 = \{ \ell_1 \}$, we set $Y = 1$.  Then~\eqref{subselimdeq} 
becomes
\begin{equation*}
1 = \sum_{\ell \in L_2 \setminus L_1} B_\ell \, (j_1 - k(\ell)) p_{j_1} q_{k(\ell)}
+ B_{\ell_1} ( (j_1-k_1) p_{j_1} q_{k_1} - c_{\ell_1}).
\end{equation*}
Substituting $q_{k_1} = 1$ and all other $q_k = 0$, leads to a polynomial $C_{\ell_1}$ obtained from $B_{\ell_1}$ by this substitution and
$1 = C_{\ell_1} ((j_1 - k_1) p_{j_1} - c_{\ell_1})$, 
which cannot 
hold for all $p_{j_1} \in \bC$, so we reached a contradiction.

If $\abs{L_2 \cap L_1} > 1$, let $\ell_2 \in (L_2 \cap L_1) \setminus \{ \ell_1 
\}$, then $p_{j_1} q_{k(\ell_2)}$ has nonzero coefficient in 
$\rho_{\ell_2}$, and we set $Y = q_{k(\ell_2)}$.
Substitute $q_k = 0$ for $k \notin \{ k(\ell) : \ell \in L_1 \cap L_2 \}$,
and denote by $D_\ell$ the polynomial obtained from $B_\ell$ by 
this substitution.
Then~\eqref{subselimdeq} becomes
\begin{equation*}
\begin{split}
q_{k(\ell_2)}
&= D_{\ell_1} ( (j_1 - k_1) p_{j_1} q_{k_1} - c_{\ell_1}) \\
&+
\sum_{\ell \in L_2 \cap L_1 \setminus \{\ell_1\}}
D_\ell \left( (j_1 - k(\ell)) p_{j_1} q_{k(\ell)} - 
\frac{c_\ell}{c_{\ell_1}} (j_1 - k_1) p_{j_1} q_{k_1} \right).
\end{split}
\end{equation*}
Finally, we substitute $p_{j_1} = 1$ and $q_{k_1} = c_{\ell_1}/(j_1 - k_1)$ and
$q_{k(\ell)} = \frac{c_{\ell}}{j_1 - k(\ell)}$ for $\ell \in L_2 \cap L_1 
\setminus \{ \ell_1, \ell_2 \}$, which,
denoting by $E_{\ell_2}$ the polynomial obtained from $D_{\ell_2}$ by the substitution,
yields
\begin{equation*}
q_{k(\ell_2)} = E_{\ell_2} ( (j_1 - k(\ell_2)) q_{k(\ell_2)} - c_{\ell_2} ).
\end{equation*}
This is impossible since $c_{\ell_2} \neq 0$ since $\ell_2 \in L_1$, and
$j_1 - k(\ell_2) \neq 0$ since $\ell_2 \in L_2$.
This concludes the proof of existence.

\emph{Step 2: Number of normalized solutions.}
By Remark~\ref{rem:critpts_poly}, the solutions can
be normalized by~\eqref{leadingcoeffassump}.
By the correspondence established in
Proposition~\ref{prop:critpts_rat_vs_poly}, the number of (normalized)
solutions in Theorems~\ref{thm:critpts_poly} and~\ref{thm:critpts}
are the same, and we show the bound for the latter.
If $F$ is a rational function and $\varphi(z) = \frac{az+b}{cz+d}$
with $ad-bc \neq 0$ a M\"obius transformation,
then $F$ and $\varphi \circ F$ have the same critical points.
Following Goldberg~\cite{Goldberg}, we call two rational functions
$F, F_1$ equivalent if $F_1 = \varphi \circ F$ with a M\"obius transformation
$\varphi$.
Goldberg~\cite[Thm.~1.3]{Goldberg} showed that the number of equivalence
classes (with respect to this equivalence relation) of rational functions
with $2n$ prescribed critical points is bounded by the Catalan number
$C_n = \frac{1}{n+1} \binom{2n}{n}$.
Given a rational function of degree $n+1$ with the critical points
$\cp_1, \ldots, \cp_{2n}$, there is an equivalent rational function of
type $(n+1, n)$ and, moreover, there is one that is normalized by
$F(z) = z + O(1/z)$ at infinity.
We show that this is the only normalized function in its equivalence class.
Indeed, if $F_1 = \varphi \circ F$ with a M\"obius transformation
$\varphi$ and $F_1(z) = z + O(1/z)$, then
$\infty = F_1(\infty) = \varphi(F(\infty)) = \varphi(\infty)$,
hence $\varphi$ has the form $\varphi(z) = az+b$ with
$a \in \C \setminus \{  0 \}$ and $b \in \C$.
Then, at infinity, $F_1(z) = z + O(1/z) = \varphi(F(z)) = a z + b + O(1/z)$,
which implies $a = 1$, $b = 0$, i.e., $\varphi(z) = z$ and $F_1 = F$.
Thus, there is only one rational function normalized by $z + 
O(1/z)$ at infinity in each equivalence class, which establishes the
bound on the number of normalized solutions in Theorem~\ref{thm:critpts}
and Theorem~\ref{thm:critpts_poly}.
\end{proof}

\begin{remark}
In~\cite[Thm.~1.3]{Goldberg}, it is also shown that the bound
in Theorems~\ref{thm:critpts} and~\ref{thm:critpts_poly}
is attained in the generic case, though not always.
For $n = 2$, the exact number is characterized in~\cite[Thm.~1.4]{Goldberg}:
There is exactly one rational function if the cross ratio of $(\cp_1, \ldots, 
\cp_4)$ is equal to $(1 + i \sqrt{3})/2$, and there are exactly $C_2 = 2$ 
such rational functions otherwise.
\end{remark}

\section{Solutions for \texorpdfstring{$n=1$}{n=1}
and symmetry considerations}
\label{sect:discussion}

If $n = 1$, the rational functions in Theorems~\ref{thm:openup},
\ref{thm:critvals} and~\ref{thm:critpts} are uniquely determined when
normalized by $F(z) = z + O(1/z)$ at infinity, and
can be determined explicitly.

\begin{proposition}
\label{prop:casen1}
\ \ \ 
\begin{enumerate}
\item \label{it:critpts_n1}
Given distinct $\cp_1, \cp_2 \in \bC$, then
$F(z) = z + \frac{r}{z - p}$ with
$r = ((\cp_1 - \cp_2)/2)^2$ and $p = (\cp_1 + \cp_2)/2$
is the unique rational function in Theorem~\ref{thm:critpts} with the 
normalization $F(z) = z + O(1/z)$ at infinity.
Its critical values are
$\cv_1 = (3 \cp_1 - \cp_2)/2$,
$\cv_2 = (3 \cp_2 - \cp_1)/2$.

\item \label{it:critvals_n1}
Given distinct $\cv_1, \cv_2 \in \bC$, then
$F(z) = z + \frac{r}{z - p}$ with
$r = ((\cv_1 - \cv_2)/4)^2$ and
$p = (\cv_1 + \cv_2)/2$
is the unique rational function in Theorem~\ref{thm:critvals} with the 
normalization $F(z) = z + O(1/z)$.
Its critical points are
$\cp_1 = (3 \cv_1 + \cv_2)/4$,
$\cp_2 = (\cv_1 + 3 \cv_2)/4$.

\item Let $\gamma$ be a Jordan arc with distinct endpoints $\eta_1, \eta_2$.
Then $F$ from \ref{it:critvals_n1}.\ is the unique open up mapping from a Jordan domain $G$ to 
$\bC_\infty \setminus \gamma$.  The boundary of $G$ consists of the two branches 
of $F^{-1}(\gamma)$.
\end{enumerate}
\end{proposition}

\begin{proof}
Let $F(z) = z + \frac{r}{z - p}$ and assume that $F(\cp_j) = \cv_j$ and 
$F'(\cp_j) = 0$ for $j=1,2$.
Then $F'(\cp_j) = 0$ is equivalent to
$r = (\cp_j - p)^2$, for $j = 1, 2$.
This implies $(\cp_1 - p)^2 = (\cp_2 - p)^2$ and $\cp_1 - p = - (\cp_2 - p)$
(because $\cp_1 - p = \cp_2 - p$ is impossible, since $\cp_1 \neq \cp_2$ in~\ref{it:critpts_n1}.,
and $\cp_1 = \cp_2$ implies $\cv_1 = \cv_2$, which is not possible in~\ref{it:critvals_n1}.).
We thus obtain $p = (\cp_1 + \cp_2)/2$ and $r = ((\cp_1 - \cp_2)/2)^2$.
To complete~\ref{it:critpts_n1}., the critical values are
$\cv_1 = f(\cp_1) = (3 \cp_1 - \cp_2)/2$ and
$\cv_2 = f(\cp_2) = (3 \cp_2 - \cp_1)/2$.
For~\ref{it:critvals_n1}., we compute $\cp_1, \cp_2$ from $\cv_1, \cv_2$, 
and obtain successively
$\cp_1 = (3 \cv_1 + \cv_2)/4$,
$\cp_2 = (\cv_1 + 3 \cv_2)/4$, 
$r = ((\cv_1 - \cv_2)/4)^2$ and
$p = (\cv_1 + \cv_2)/2$.
Finally, if $\gamma$ is a Jordan arc with endpoints $\cv_1, \cv_2$, then $F$
in~\ref{it:critvals_n1}.\ 
is the open up mapping, and the domain $G$ is bounded by the two branches of
$F^{-1}(\gamma)$.
\end{proof}

An open up mapping is also a solution of the prescribed critical
values problem; see Proposition~\ref{prop:Prob1vsProb2}.
Conversely, a solution of the critical value problem is in general not an open 
up mapping when $n \geq 2$, and we give 
two examples
below.
The reason becomes apparent from Theorem~\ref{thm:deformation_of_arcs}:
The open up mapping in Theorem~\ref{thm:openup} 
depends on the endpoints of the arcs 
(critical values) and the topology of $\bC \setminus \cup_{j=1}^n \gamma_j$.
In contrast, only the information about the critical values is present in 
Theorem~\ref{thm:critvals}.
For $n = 1$, 
a rational function of type $(2, 1)$ has critical values $\cv_1, \cv_2$
if and only if it is an open up mapping for any Jordan arc connecting
$\cv_1, \cv_2$.

We reformulated Theorems~\ref{thm:critvals} and~\ref{thm:critpts}
in terms of polynomials.
Since the number of normalized solutions is finite, 
one can apply several existing methods and algorithms to
solve the polynomial systems symbolically.
We just mention~\cite[Sect.~10.4]{CLOIdeals}, \cite{EderFaugere}
and the references therein.
For many of the methods, it is crucial that the number of solutions is finite.
Also, there are several numerical algorithms and solvers available, see
e.g.~\cite{AllgowerGeorg} and~\cite{SommeseWampler2005}.
Finally, if the finitely many solutions of the critical
value problem have been obtained, one can go through these
solutions and verify which one is the open up mapping
for a given set of arcs.
We use this approach in Section~\ref{sect:example}.

\medskip

The examples in the following section show that it is advantageous to 
simplify the polynomial equations, e.g. by reducing the number 
of equations and unknowns.  Therefore we investigate some 
symmetric settings and how one can simplify the polynomial equations.
We collect three useful results on the open up mapping when $E$ has some symmetry.

\begin{lemma} \label{lem:openupmap_odd}
Let $E = \gamma_1 \cup \ldots \cup \gamma_n$ be the union of the disjoint 
Jordan arcs $\gamma_1, \ldots, \gamma_n$, and let $F$ be the
open up mapping of type $(n+1, n)$ with $F(z) = z + O(1/z)$ at infinity.
\begin{enumerate}
\item Suppose that $z \in E$ if and only if $-z \in E$. 
Then $F$ is odd, i.e.,~$F(-z)=-F(z)$.
\item \label{item:openupmap_odd_2}
Suppose that $z \in E$ if and only if $\overline{z} \in E$.
Then $F$ is real, i.e., $F(z) = \overline{F(\overline{z})}$.
\end{enumerate}
\end{lemma}

\begin{proof}
Let $F : \bC_\infty \setminus K \to \bC_\infty \setminus E$ with $F(z) = z + 
O(1/z)$ be the open up mapping from Theorem~\ref{thm:openup}.
Then
\begin{equation*}
G : \bC_\infty \setminus (-K) \to \bC_\infty \setminus E, \quad G(z) = - F(-z),
\end{equation*}
is also an open up mapping of type $(n+1, n)$ with 
$G(z) = z + O(1/z)$ for $z \to \infty$, hence $F = G$ by 
Theorem~\ref{thm:openup}.
The proof of~\ref{item:openupmap_odd_2} is similar.
\end{proof}

\begin{lemma} \label{lem:odd_F}
Let $F(z) = z + \sum_{j=1}^n \frac{r_j}{z - p_j}$ with 
$r_1, \ldots, r_n \in \bC \setminus \{ 0 \}$ and distinct poles 
$p_1, \ldots$, $p_n \in \bC$.
Then $F$ is odd if and only if the poles appear in pairs $\pm p$ with equal 
residues.
In particular, if $F$ and $n$ are odd, then one pole is at the origin.
\end{lemma}

\begin{proof}
Let $F$ be odd. Since the partial fraction decomposition is unique and
\begin{equation*}
F(z) = - F(-z) = z + \sum_{j=1}^n \frac{r_j}{z+p_j},
\end{equation*}
the poles appear in pairs $\pm p$ with equal residues.
The converse is obvious.
\end{proof}

\begin{corollary} \label{cor:crit_val_odd_solution}
Let $\cv_1, \ldots, \cv_{2n} \in \bC$ be distinct such that $\cv \in
\{ \cv_1, \ldots, \cv_{2n} \}$ implies $-\cv \in \{ \cv_1, \ldots,
\cv_{2n} \}$.  Then there is an odd rational function that solves
the critical value problem in Theorem~\ref{thm:critvals}.
\end{corollary}

\begin{proof}
Connect the critical values $\cv_1, \ldots, \cv_{2n}$ by arcs $\gamma_1,
\ldots, \gamma_n$, such that $E = \cup_{j=1}^n \gamma_j$ has the property
$z \in E$ if and only if $-z \in E$.
Note that the arcs themselves need not be symmetric with respect
to the origin, but their union $E$ needs to be symmetric.
The arcs can be constructed as follows.  
There exists a line through the origin that divides the plane
in two half-planes, each containing $n$ of the points $\cv_1, \ldots, \cv_{2n}$ (by symmetry).
If $n \geq 2$, choose two points in one half-plane and connect
them by a Jordan arc $\gamma_1$ in that half-plane, such that
$\gamma_1$ contains no other critical value.
Then $\gamma_2 = - \gamma_1$ also connects two critical values
in the other half-plane.  This way, we construct disjoint Jordan
arcs connecting the critical values.
If there is only one critical point left in each half-plane (which happens
if $n$ is odd),
we connect these two by a Jordan arc $\gamma_n$ that is symmetric
with respect to the origin (and disjoint from all previous arcs).
Then the open up mapping of $E$ is odd by
Lemma~\ref{lem:openupmap_odd} and a solution of the critical
value problem.
\end{proof}

\section{Two examples and further comments}
\label{sect:example}

We give two examples when $n=2$, i.e., in the case of two
disjoint Jordan arcs $\gamma_1, \gamma_2$ with endpoints
$\cv_1, \cv_2, \cv_3, \cv_4$.
The open up mapping in this case has the form
\begin{equation} \label{eqn:general_F}
F(z) = z + \frac{r_1}{z - p_1} + \frac{r_2}{z - p_2}
\end{equation}
with $r_1, r_2 \in \bC \setminus \{ 0 \}$ and distinct $p_1, p_2 \in \bC$.
Without loss of generality, we assume here and in the following that the
rational functions
are normalized by $z + O(1/z)$ at infinity; see Theorems~\ref{thm:openup}
and~\ref{thm:critvals} and also Remark~\ref{rem:lineartransformation}.

We also determine
all rational functions of the form~\eqref{eqn:general_F} with critical values $\cv_1, \cv_2, \cv_3, \cv_4$.  
The problem then is to find $r_1, r_2, p_1, p_2$ and 
$\cp_1, \cp_2, \cp_3, \cp_4 \in \bC$ with
\begin{equation}
F(\cp_j) = \cv_j, \quad F'(\cp_j) = 0, \quad j = 1, 2, 3, 4.
\label{eq:simpleexamplefortwo}
\end{equation}
By Theorem~\ref{thm:critvals}, there are $3 H_2 = 12$ rational functions
of the form~\eqref{eqn:general_F} satisfying~\eqref{eq:simpleexamplefortwo}.
By Proposition~\ref{prop:Prob1vsProb2}, the open up mapping is among these.

\subsection{An example with symmetry}

As first example, we determine the open up mapping for
\begin{equation} \label{eqn:symmetric_E}
E = \cc{-2, -1} \cup \cc{1, 2}
\end{equation}
with Jordan arcs $\gamma_1 = \cc{1, 2}$ and $\gamma_2 = \cc{-2, -1}$
as well as all rational functions of the form~\eqref{eqn:general_F}
with critical values
\begin{equation}
\cv_1 = 1, \quad \cv_2 = 2, \quad \cv_3 = -1, \quad \cv_4 = -2.
\label{eq:symmetricexamplevalues}
\end{equation}
Since $E$ is symmetric with respect to the origin,
the open up mapping~\eqref{eqn:general_F} has the simpler form
\begin{equation} \label{eqn:F_odd_2poles}
F(z) = z + \frac{r}{z - p} + \frac{r}{z + p}
\quad \text{with } r, p \in \bC \setminus \{ 0 \}
\end{equation}
by Lemma~\ref{lem:openupmap_odd} and Lemma~\ref{lem:odd_F}.

\begin{figure}[t]
{\centering
\includegraphics[width=0.45\linewidth]{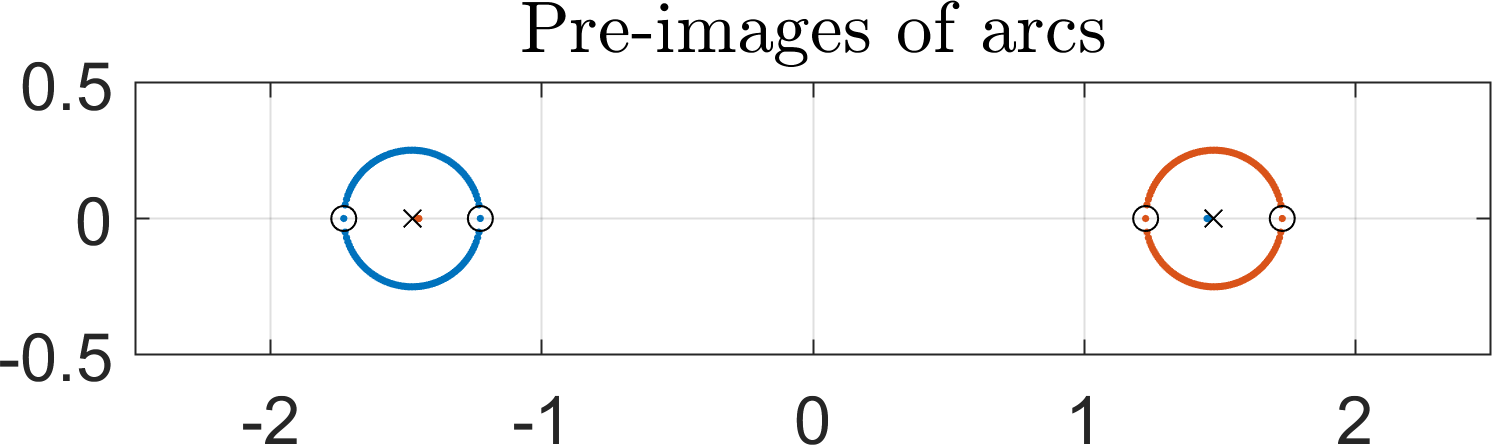}
\includegraphics[width=0.45\linewidth]{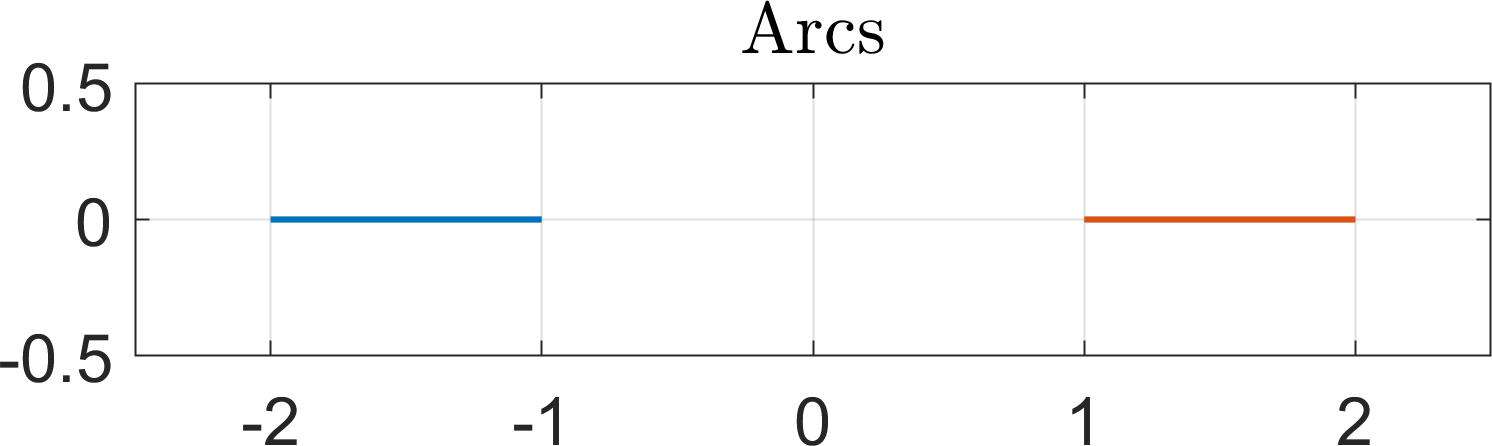}

}
\caption{The function in~\eqref{eqn:ex_openup} is the minimal degree rational open up 
mapping for $\cc{-2, -1} \cup \cc{1, 2}$.  Right: Arcs.  Left: Pre-images of 
the arcs, critical points (black circles) and poles (crosses) of $F_2$.}
\label{fig:ex_openup}
\end{figure}

Accordingly, we simplify also the critical value problem 
and ask only for odd solutions $F$ of the form~\eqref{eqn:F_odd_2poles}
of
\begin{equation} \label{eqn:sym_crit_val_pblm}
F(\pm \cp_1) = \pm \cv_1 = \pm 1, \quad F(\pm \cp_2) = \pm \cv_2 = \pm 2,
\quad F'(\pm \cp_1) = F'(\pm \cp_2) = 0.
\end{equation}
One can obtain
four distinct symbolical solutions of~\eqref{eqn:sym_crit_val_pblm}
(the following numbers are rounded to $4$ digits)
\begin{align}
F_1(z) &= z + \frac{0.1272}{z - 0.0658i} + \frac{0.1272}{z + 0.0658i}, \\
F_2(z) &= z + \frac{0.0630}{z - 1.4786} + \frac{0.0630}{z + 1.4786}, 
\label{eqn:ex_openup} \\
F_3(z) &= z + \frac{0.4605 + 0.1279i}{z - (0.3958 - 0.3693i)}
+ \frac{0.4605 + 0.1279i}{z + (0.3958 - 0.3693i)}, \\
F_4(z) &= z + \frac{0.4605 - 0.1279i}{z - (0.3958 + 0.3693i)}
+ \frac{0.4605 - 0.1279i}{z + (0.3958 + 0.3693i)}.
\end{align}
By Lemma~\ref{lem:openupmap_odd}, the open up map of $E = \cc{-2, -1} \cup \cc{1, 2}$ is odd and real (i.e., $F(z) = \overline{F(\overline{z})}$),
and hence must be $F_2$.
By computing the pre-images of $\gamma_1$ and $\gamma_2$, we find
that $F_2$ is indeed
the (unique) open up mapping for $\gamma_1$ and $\gamma_2$,
see Figure~\ref{fig:ex_openup}, while the other 
functions are not open up mappings for $\gamma_1$ and $\gamma_2$.
This confirms that a function with critical values at the endpoints of the arcs
in Theorem~\ref{thm:critvals} is in general
not an open up mapping for the arcs.

\begin{figure}[t]
{\centering
\includegraphics[width=0.30\linewidth]{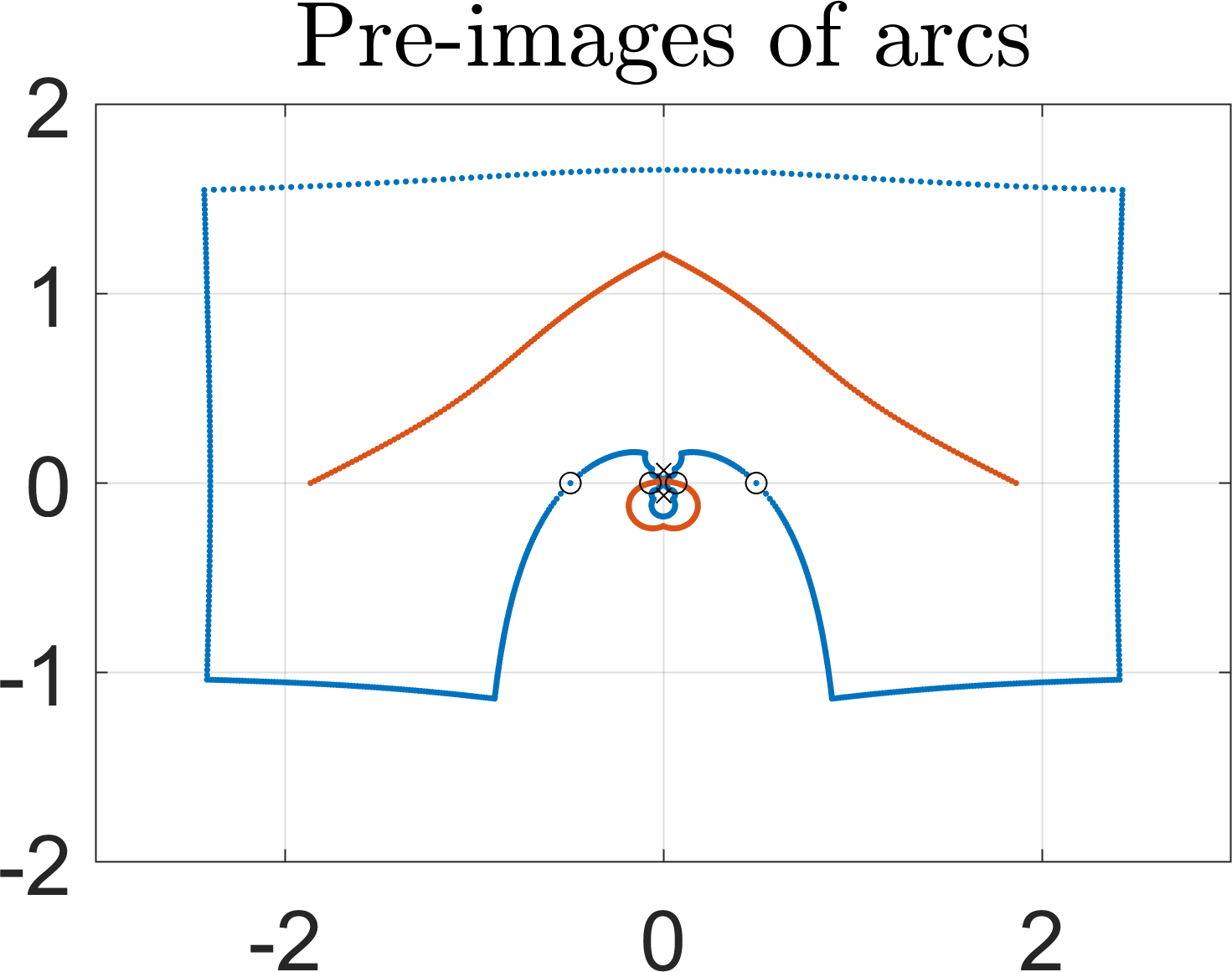}
\includegraphics[width=0.30\linewidth]{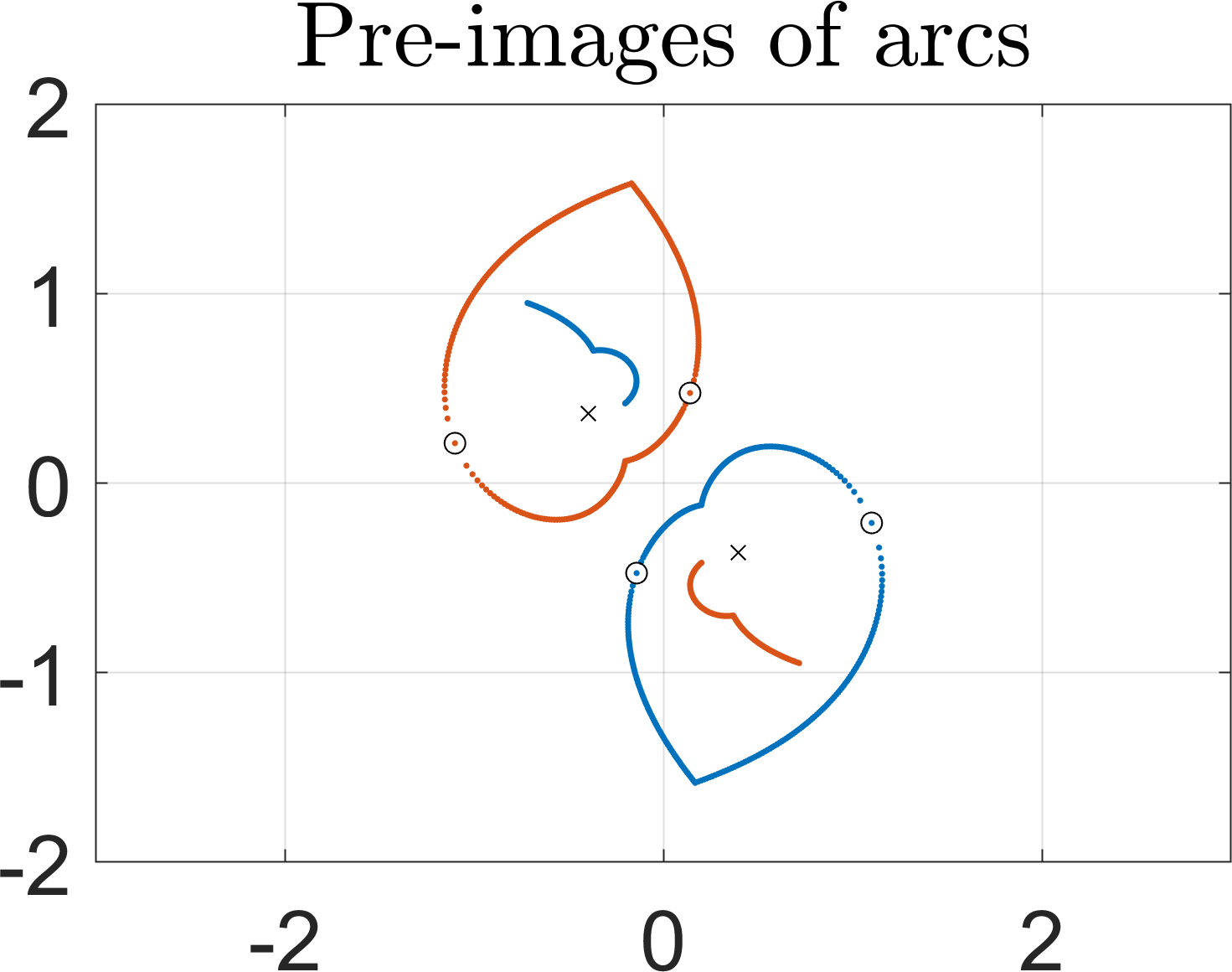}
\includegraphics[width=0.30\linewidth]{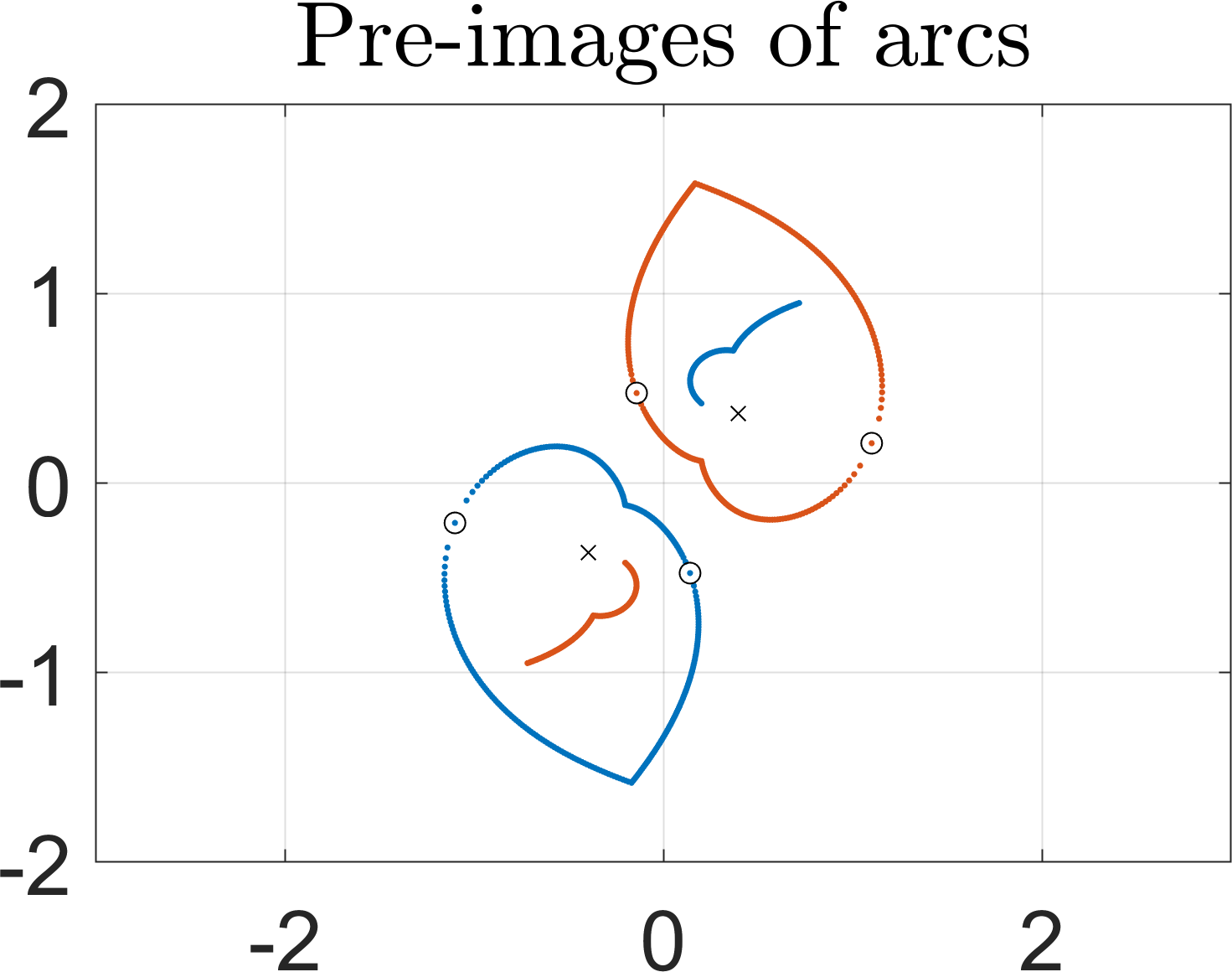}

\includegraphics[width=0.30\linewidth]{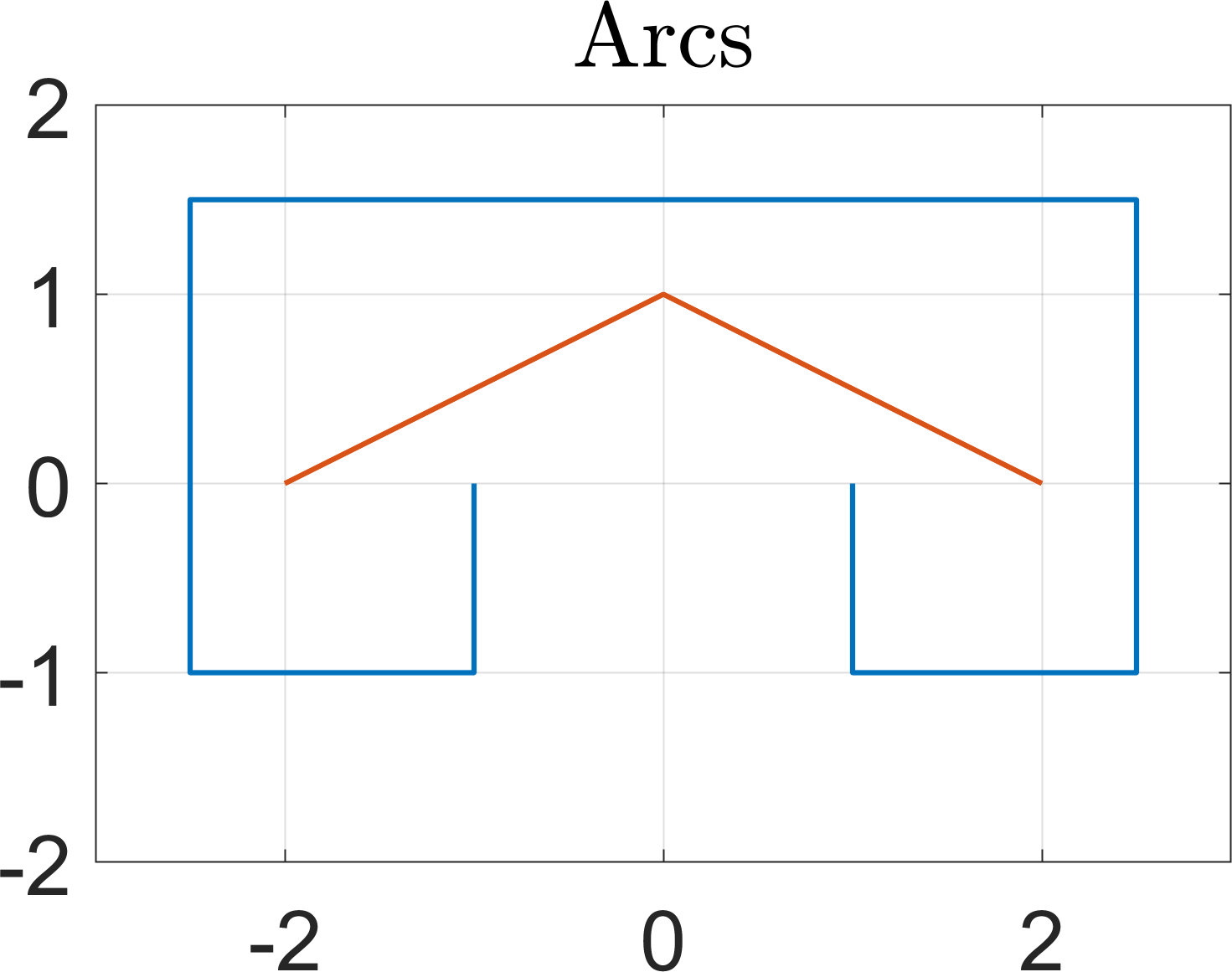}
\includegraphics[width=0.30\linewidth]{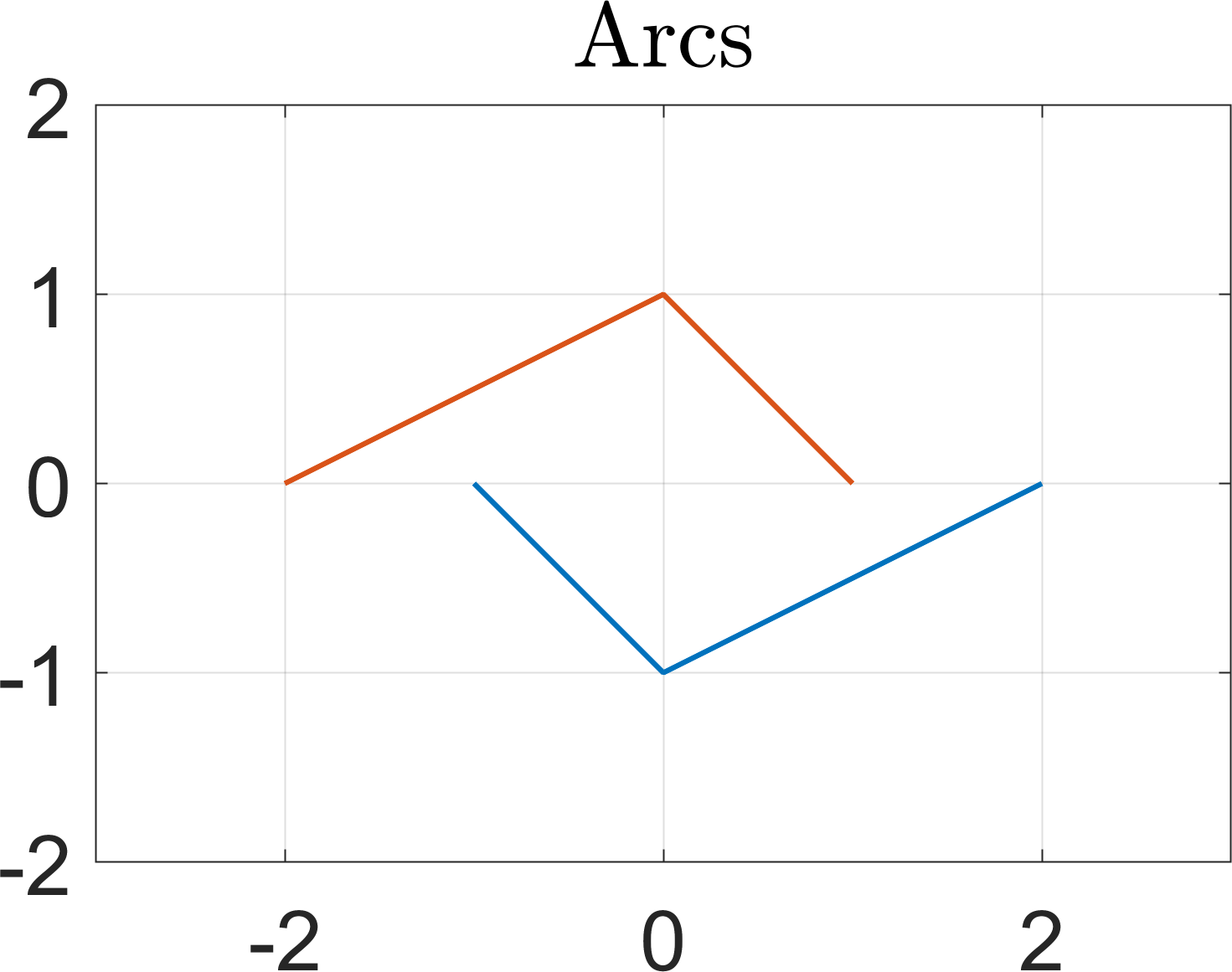}
\includegraphics[width=0.30\linewidth]{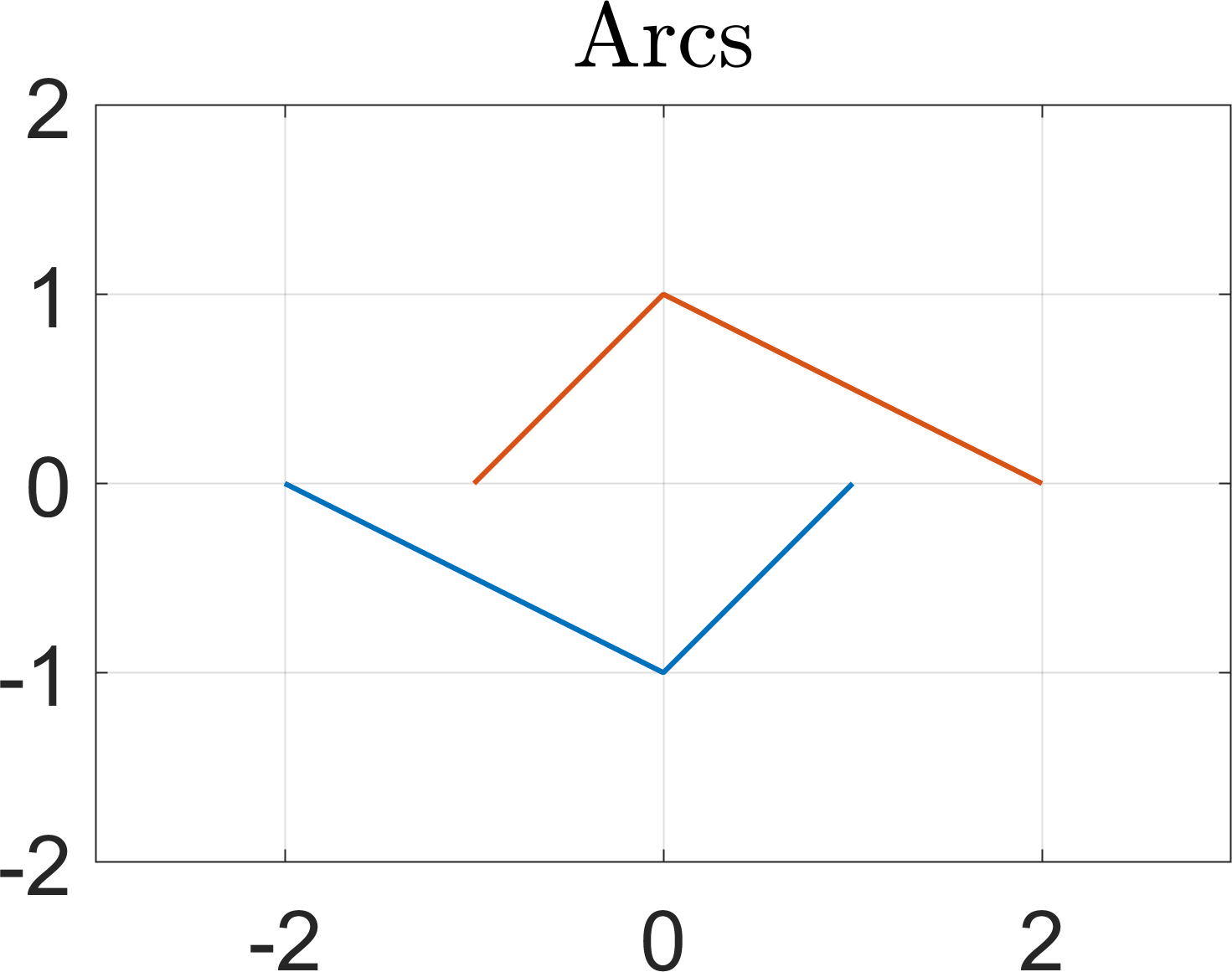}

}
\caption{Arcs (bottom) and pre-images (top) that are opened up under 
$F_1, F_3, F_4$ (left to right).}
\label{fig:arcs_and_preimages}
\end{figure}

However, each of the functions $F_1, F_2, F_3, F_4$ is an open up mapping for
a suitable set of 
disjoint Jordan arcs, where each Jordan arc connects two of the critical 
values.  Examples of configurations opened up by $F_1, F_3, F_4$ are displayed 
in Figure~\ref{fig:arcs_and_preimages}.
This leads to the following conjecture.

\begin{conjecture}
Let $F$ be a rational function of type $(n+1,n)$ with distinct critical
values $\cv_1, \ldots, \cv_{2n} \in \bC$.
Then there exists a set of disjoint Jordan arcs $\gamma_1,\ldots,\gamma_n$,
each arc connecting two points in $\{ \cv_1, \ldots, \cv_{2n} \}$,
such that $F$ is the open up mapping in Theorem~\ref{thm:openup}
for the arcs $\gamma_1,\ldots,\gamma_n$.
\end{conjecture}

A deformation (fixed endpoint homotopy) of the arcs as described in 
Theorem~\ref{thm:deformation_of_arcs} yields homotopic configurations that are 
also opened up by the same function.
Note that for some functions there are also other (non-homotopic) 
configurations that are opened up.  This shows that the open up configuration 
for a given function in Theorem~\ref{thm:critvals} is not unique up to 
fixed endpoint homotopy.

Next, we consider the prescribed critical value problem
in general form, i.e.,~\eqref{eqn:general_F}
and~\eqref{eq:simpleexamplefortwo}.
We multiply~\eqref{eq:simpleexamplefortwo} with the denominators
and solve the obtained polynomial equations symbolically (with Singular and also with Magma),
The result of this computation is the uniquely determined Gröbner basis 
consisting of $41$ polynomials.
This new set of polynomial equations is then solved.
Some of the coefficients (which are rational numbers) have 
numerators and denominators in reduced form of magnitude $10^{40}$.

The computation yields $12$ distinct rational functions,
i.e., all solutions,
including the odd functions $F_1, \ldots$, $F_4$ above, and $8$ solutions
that are not of the form~\eqref{eqn:F_odd_2poles}.
In particular, the symmetry of the critical values (endpoints of the arcs)
does not imply symmetry of the rational function.
The symmetric problem of prescribed 
critical values (with values $\pm 1$ and $\pm 2$, see \eqref{eqn:general_F} 
and \eqref{eq:symmetricexamplevalues}) 
has a nonsymmetric solution 
(e.g., 
$r_1=-0.0005$, $r_2=0.0630$, $p_1=-1.4786$, $p_2=-1.4998$
where the values are rounded to $4$ digits).

In this first example, we computed the open up mapping of
$\cc{-2, -1} \cup \cc{1, 2}$ by solving the simplified prescribed critical
value problem~\eqref{eqn:sym_crit_val_pblm},
and also computed all rational functions with prescribed
critical values $1, 2, -1, -2$.
If the aim is to compute the open up mapping, then
\eqref{eqn:sym_crit_val_pblm} is better suited than~\eqref{eq:simpleexamplefortwo},
since it involves approximately half the number of unknowns and equations
and is thus easier and faster to solve.

\subsection{An example without symmetry}

In the second example let 
$\gamma_1=[1,2]$ and
$\gamma_2=[i,2i]=\{(1-t)i+t 2i:\ 0\le t\le 1\}$.
Therefore we consider
the critical value 
problem~\eqref{eq:simpleexamplefortwo} with
\begin{equation} \label{eqn:crit_vals_ex2}
\eta_1=1, \quad  \eta_2=2, \quad
\eta_3=i, \quad  \eta_4=2i.
\end{equation}

We determined the corresponding Gröbner basis with Singular 
and Magma.
The Gröbner basis consists of
$77$ polynomials and some of the rational coefficients have numerators and 
denominators of size $10^{75}$.

\begin{figure}
{\centering
\includegraphics[width=0.32\linewidth]{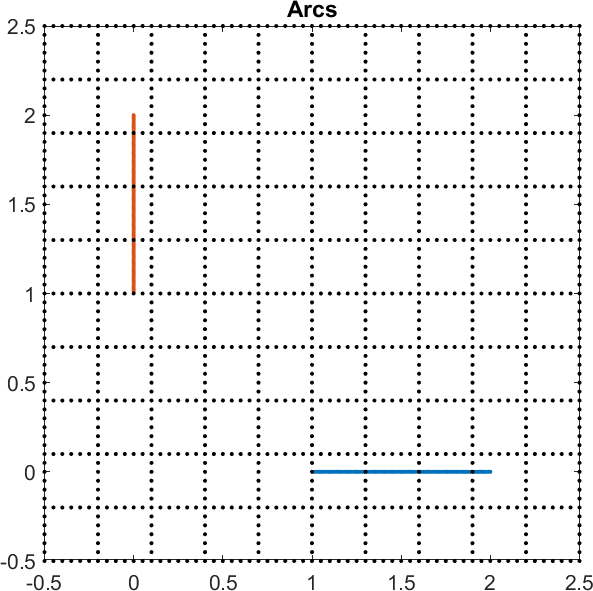}
\includegraphics[width=0.32\linewidth]{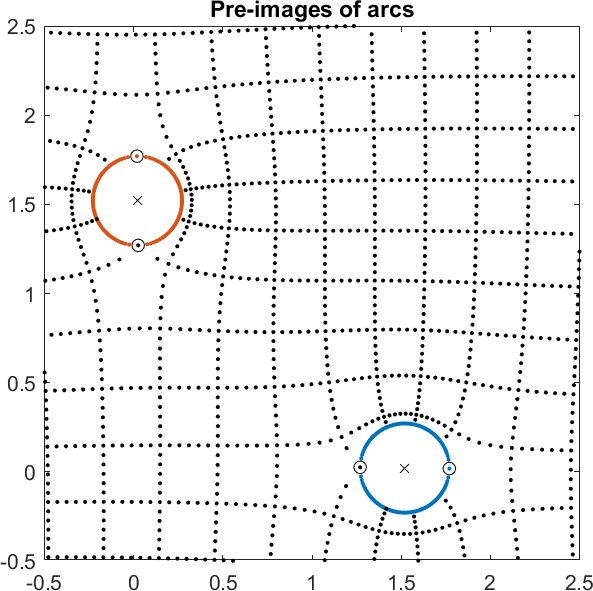}

}
\caption{Open up mapping for $\cc{1, 2} \cup \cc{i, 2i}$.
Left: arcs and a grid (black).
Right: Domain $\bC_\infty \setminus K$ bounded by Jordan curves and the 
pre-image of the grid under the open up mapping.}
\label{fig:nonsymmetric_example}
\end{figure}

In total, we obtain $12$ distinct rational functions of the
form~\eqref{eqn:general_F} with critical values $1, 2, i, 2i$,
i.e., all solutions have been computed; see Theorem~\ref{thm:critvals}.
Finally, we determine which solution is the open up mapping for
$\cc{1, 2} \cup \cc{i, 2i}$ by computing the pre-images of
$\gamma_1 = \cc{1, 2}$ and $\gamma_2 = \cc{i, 2i}$, and obtain that
the open up mapping is (coefficients rounded to four digits)
\begin{equation}
F(z) = z + \frac{-0.0625 - 0.0009i}{z - (0.0214 + 1.5203i)} + \frac{0.0625 - 0.0009i}{z - (1.5203 + 0.0214i)}.
\end{equation}
Figure~\ref{fig:nonsymmetric_example} visualizes the open up mapping.
The left panel shows the arcs (blue and red) with a grid (black dots).
The right panel shows the pre-image domain $\bC_\infty \setminus K$
bounded by Jordan curves (blue and red) and the pre-image of the grid
under $F$ (black dots).

\section*{Funding}
This work was supported by Moscow Center for Fundamental and Applied Mathematics, Agreement with the Ministry of Science and Higher Education of the Russian Federation [grant no.~075-15-2022-283] to S.K.

This work was supported in part by the János Bolyai 
Scholarship of the Hungarian Academy of Sciences to B.N.

This research was supported by project TKP2021-NVA-09. Project no.~TKP2021-NVA-09 has been implemented with the support provided by the Ministry of Innovation and Technology of Hungary from the National Research, Development and Innovation Fund, financed under the TKP2021-NVA funding scheme.

\section*{Acknowledgments}

Béla Nagy would like to thank Mikhail Tyaglov and Sergei Kalmykov from SJTU 
for the hospitality and support during visits to Shanghai.

We are grateful to Vilmos Totik for valuable comments which helped to improve the
presentation of the paper.
We thank Vladimir Lysov for helpful discussions and suggesting the bound on
the number of solutions in Theorem~\ref{thm:critvals}.

We thank the anonymous referee for helpful and constructive comments that
lead to improvements in the paper.

\bibliographystyle{amsplain}
\bibliography{amsrefs}
\medskip

Sergei Kalmykov\\
School of Mathematical Sciences,\\
CMA-Shanghai,\\
Shanghai Jiao Tong University\\
800 Dongchuan RD\\
Shanghai, 200240\\
P.R. China \\
\href{mailto:kalmykovsergei@sjtu.edu.cn}{kalmykovsergei@sjtu.edu.cn}
\\and
\\
Keldysh Institute of Applied Mathematics \\
of Russian Academy of Sciences
\\
Miusskaya pl., 4,
125047, Moscow, \\
Russia

\medskip

Béla Nagy \\
Department of Analysis,\\
Bolyai Institute\\
University of Szeged\\
Aradi v. tere 1\\
6720 Szeged \\
Hungary  \\
\href{mailto:nbela@math.u-szeged.hu}{nbela@math.u-szeged.hu}

\medskip

Olivier S\`ete \\
Institute of Mathematics and Computer Science \\
Universit\"at Greifswald \\
Walther-Rathenau-Stra{\ss}e 47 \\
17489 Greifswald \\
Germany \\
\href{mailto:olivier.sete@uni-greifswald.de}{olivier.sete@uni-greifswald.de}

\end{document}